\documentclass[12pt]{iopart}

\usepackage{paralist}
\usepackage{graphics} 
\usepackage{epsfig} 
\usepackage{graphicx}  \usepackage{epstopdf}
\usepackage{mathrsfs}

\usepackage[breaklinks,colorlinks=true,linkcolor=blue,citecolor=red,backref=page]{hyperref}

\newtheorem{thm}{Theorem}[section]
\newtheorem{cor}[thm]{Corollary}

\newtheorem{proposition}[thm]{Proposition}
\newtheorem{example}[thm]{Example}

\newtheorem{remark}[thm]{Remark}

\newcommand{\eqref}[1]{(\ref {#1})}

\newcommand{\ea}{\end{eqnarray}}  
\newcommand{\ba}{\begin{eqnarray}}  
\newcommand{\ean}{\end{eqnarray*}}  
\newcommand{\ban}{\begin{eqnarray*}}

\DeclareMathAlphabet{\itbf}{OML}{cmm}{b}{it}
\DeclareMathAlphabet\mathbfcal{OMS}{cmsy}{b}{n}

\def\qed{\hfill {\small $\Box$} \\}  
 
\def\EE{\mathbb{E}}
 
 \newcommand{\dessous}[2]{
\renewcommand{\arraystretch}{0.5} 
\begin{array}[t]{c}
{#1} \\
\scriptstyle
{#2}
\displaystyle
\end{array}
\renewcommand{\arraystretch}{1.0}
}

\usepackage{iopams}

\begin{document}

\title{Low-frequency source imaging in an acoustic waveguide}
\author{Josselin Garnier$\hbox{}^{(1)}$}
\address{$\hbox{}^{(1)}$ 
CMAP, CNRS, Ecole Polytechnique, Institut Polytechnique de Paris, 
  Palaiseau, France.}
  \ead{josselin.garnier@polytechnique.edu}
 \vspace{10pt}
\begin{indented}
\item[] 
\end{indented}

\begin{abstract}
Time-harmonic far-field source array imaging in a two-dimensional waveguide is analyzed.
A low-frequency situation is considered in which the diameter of the waveguide is slightly larger than the wavelength,
so that the waveguide supports a limited number of guided modes,
and the diameter of the antenna array is smaller than the wavelength, so that the standard resolution formulas in open media
predict very poor imaging resolution.
A general framework to analyze the resolution and stability performances of such antenna arrays is introduced.
It is shown that planar antenna arrays perform better (in terms of resolution and stability
with respect to measurement noise) than linear (horizontal or vertical) arrays
and that vertical linear arrays perform better than horizontal arrays, for a given diameter.
However a fundamental limitation to imaging in waveguides is identified that is due to the form of the dispersion relation. 
It is intrinsic to scalar waves, whatever the complexity of the medium and the array geometry.
\end{abstract}

\vspace{2pc}

\maketitle

\section{Introduction}
We present a theoretical and numerical study of source imaging  in two-dimensional waveguides, 
using an array of sensors that record acoustic waves.  
Source imaging in waveguides is of particular interest in underwater acoustics \cite{buchanan,jensen11}. 
In a closed waveguide the wavefield can be decomposed into a finite number of guided modes and an infinite number of evanescent modes. 
In an open waveguide the wavefield can be decomposed into a finite number of guided modes and an infinite number of radiating and evanescent modes. 
The evanescent, resp. radiating, mode components of the wavefield are in general vanishing and not usable in the measured far-field data because they decay exponentially,
resp. algebraically, with the propagation distance. 
The guided mode amplitudes can be extracted from the measured data if the array is large enough
and one can then propose an imaging method that exploits them. 
The idea of formulating the inverse problem in terms of the guided mode amplitudes has recently been considered by several authors,
for source imaging \cite{bgt15} and for scatterer imaging \cite{arens11,bourgeois08,dediu06,monk12,pincon}.
However, the extraction of the guided mode amplitudes becomes challenging when the array is small \cite{tsogka13,tsogka16,tsogka18}.

In underwater acoustics, it is possible to deploy an antenna array in the oceanic waveguide but
the aperture of the array is usually limited. 
This issue is  critical when addressing low-frequency signals whose wavelengths are 
of the same order as the diameter of the waveguide so that 1) there is only a small number of guided modes
and 2) the array diameter is smaller than the wavelength.
This is typically the configuration we have in mind in this paper.
We introduce a general framework to analyze the performances (in terms of resolution and stability) of such antenna arrays.
Under ideal circumstances (i.e. in the absence of noise) the data collected by
an antenna array covering a limited part of the cross section of a waveguide can be manipulated and
processed to transform them into the set of data that would have been collected by a vertical antenna 
covering the full cross section of the waveguide, which gives full access to the guided mode amplitudes.
We explain this processing in detail in this paper.
In more realistic configurations (i.e. in the presence of noise) the processing can become unstable and requires appropriate regularization,
the imaging performance is determined by the effective rank of an operator, which depends on the array geometry and the noise level, and we analyze
different types of antennas. We show that, for a given diameter, 
planar antenna arrays perform much better (in terms of stability
with respect to measurement noise) than linear (vertical or horizontal) arrays,
and that vertical linear arrays perform better than horizontal linear arrays.
However we exhibit and clarify a fundamental limitation to imaging  in waveguides that is due to the form of the dispersion relation
and that is intrinsic to scalar waves, whatever the complexity of the medium and the array geometry.

The paper is organized as follows. In Sections \ref{sec:geo} and \ref{sec:imag} we describe the waveguide geometry 
and source array imaging problem.
In Section \ref{sec:estim} we show how to estimate the guided mode amplitudes from the array data.
In Section \ref{sec:large} we address the case of large and dense antenna arrays (large means larger than the wavelength).
In Section \ref{sec:small} we address in detail the case of small and discrete antenna arrays and consider different array geometries.

\begin{figure}
\begin{center}
\includegraphics[width=6.cm]{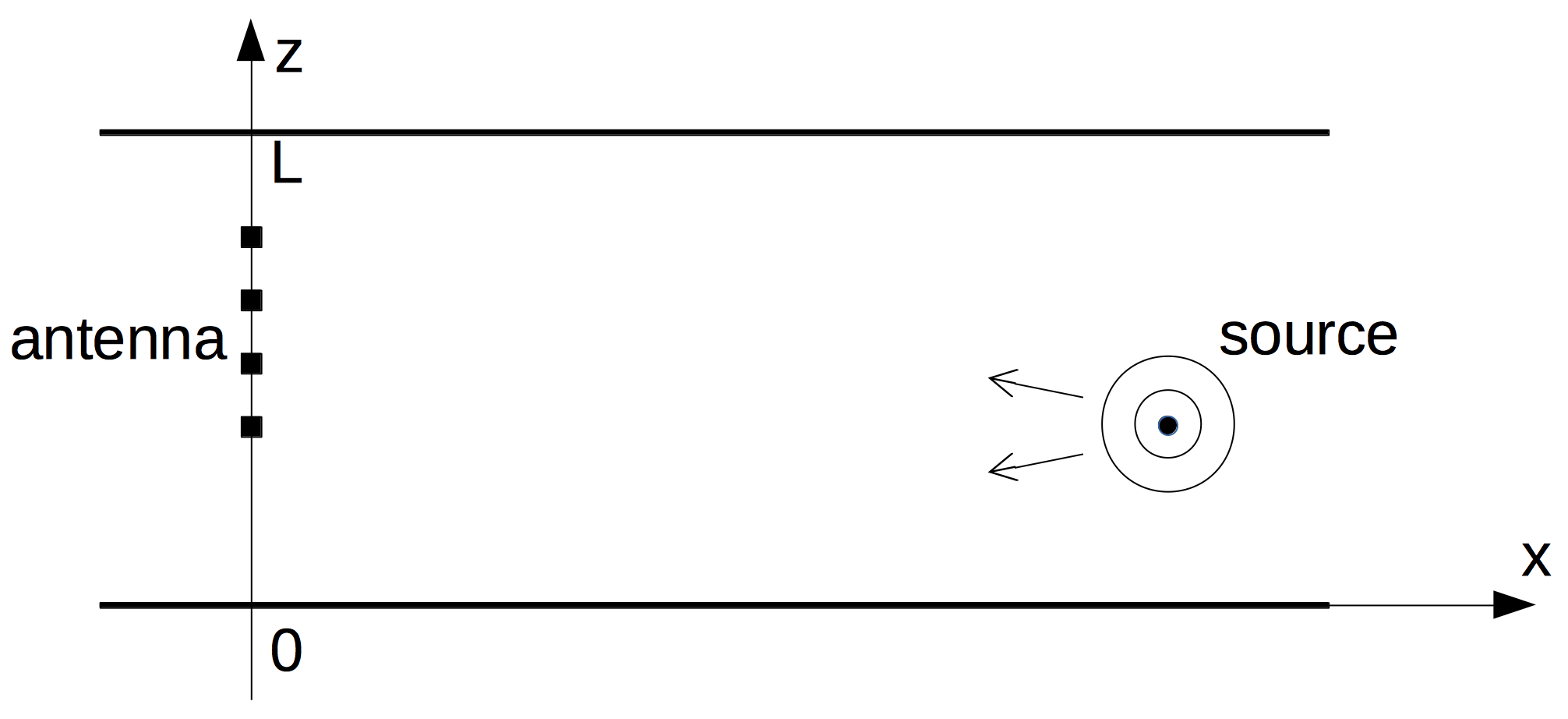} 
\end{center}
\caption{Schematic of the waveguide configuration.}
\label{fig:scheme}
\end{figure}

\section{Waveguide geometry}
\label{sec:geo}
Let us consider a two-dimensional waveguide, 
whose axis is along the $x$-axis, and the cross section is $z\in [0,L]$ (see figure \ref{fig:scheme}).
For the sake of simplicity we may consider a Dirichlet condition  at  $z=L$ (free surface)
and a Neumann or Dirichlet condition at $z=0$ (bottom).
The forthcoming results can be extended to arbitrary closed or open waveguides, such as Pekeris waveguides.
The index of refraction can be constant or variable but it depends only on $z$.
The wavefield transmitted by a time-harmonic source $s(x,z)$ at frequency $\omega$ satisfies the Helmholtz equation
\begin{equation}
\Big( \frac{d^2}{dx^2} +\frac{d^2}{dz^2} \Big)p(x,z) + \frac{\omega^2}{c^2(z)} p(x,z) = - s(x,z) , \quad (x,z) \in \mathbb{R} \times (0,L)
\end{equation}
subjected to the appropriate boundary conditions at $z=0,L$.

The eigenmodes (real-valued and orthonormal) and eigenvalues (real-valued)
of the self-adjoint operator $-\partial_z^2 - \omega^2/c^2(z)$ at frequency $\omega$ are denoted by $\phi_j(z)$ and $-\lambda_j$:
\begin{equation}
\label{eq:vp}
\frac{d^2}{dz^2}\phi_j(z) + \frac{\omega^2}{c^2(z)} \phi_j (z) =  \lambda_j \phi_j(z)  .
\end{equation}
There are $N$ guided modes for which $\lambda_j>0$ and we set $\beta_j=\sqrt{\lambda_j}$, $j=1,\ldots,N$.
The other modes for which $\lambda_j<0$  are evanescent (i.e. their amplitudes decay exponentially in $x$).
We assume thoughout the paper that the frequency $\omega$ is such that $N\geq 1$.

\section{Source imaging}
\label{sec:imag}
We consider the case of an antenna array localized in the neighborhood of the plane $x=0$.
We assume that the antenna array is supported in the domain ${\cal A} \subset [-a,a]\times [0,L]$.
The domain ${\cal A}$ can be:\\
(i) a finite collection of points $\{(x_k,z_k),k=1,\ldots,M\}$ (discrete array),\\
(ii) a square $[-a,a] \times [z_{\rm a}-a,z_{\rm a}+a]$ (continuum approximation of a dense planar array),\\
(iii) a vertical line $\{0\}\times [z_{\rm a}-a,z_{\rm a}+a]$ (continuum approximation of a dense linear vertical array localized at depth $z_{\rm a}$),\\
(iv) a horizontal line $[-a,a]\times \{z_{\rm a}\}$ (continuum approximation of a dense linear horizontal array  localized at depth $z_{\rm a}$).\\
We present a unified approach of these cases and we remark that this approach can be readily extended to other cases.
In each case we associate a corresponding uniform measure $\mu(d{\itbf x})$ with unit mass over ${\cal A}$,
such that for any test function $f$:
\begin{equation}
\int_{\cal A} f(x,z) \mu(d{\itbf x}) := 
\left\{
\begin{array}{ll}
\displaystyle
\frac{1}{M} \sum_{k=1}^M  f(x_k,z_k) ,&\mbox{ case (i),}\\
\displaystyle
\frac{1}{4a^2}  \int_{[-a,a]^2} f(x,z_{\rm a}+z) dx dz,&\mbox{ case (ii),}\\
\displaystyle
\frac{1}{2a}  \int_{-a}^{a} f(0,z_{\rm a}+z) dz,&\mbox{ case (iii),}\\
\displaystyle
\frac{1}{2a}  \int_{-a}^{a} f(x,z_{\rm a}) dx,&\mbox{ case (iv).} 
\end{array}
\right.
\end{equation}

A time-harmonic acoustic signal is transmitted by a distant source localized in the region $x>a$ (see figure \ref{fig:scheme}). 
The recorded signal is 
\begin{equation}
\label{eq:expressu0}
p(x,z)= 
\sum_{j=1}^N a_{j,o} \phi_j(z) \exp(-i \beta_j x),  \quad (x,z) \in {\cal A},
\end{equation}
where the mode amplitudes $a_{j,o}$ are determined by the source and where 
we have not written the evanescent modes, which is justified when the distance from the source to the antenna array is much larger than the wavelength.
This expression shows that the maximal information about the source 
available in the data $(p(x,z))_{(x,z)\in {\cal A}}$ 
recorded by the antenna array is the vector ${\itbf a}_o=(a_{j,o})_{j=1}^N$. The imaging procedure can be decomposed
into two steps: 1) estimation of the vector ${\itbf a}_o$ and 2) exploitation of the estimated vector to localize the source.

If we can obtain an estimate ${\itbf a}= (a _j)_{j=1}^N$ of the vector ${\itbf a}_o$ from the data $(p(x,z))_{(x,z)\in {\cal A}}$ , then
we can migrate the vector ${\itbf a}$ in order to localize the source in the region $x>0$ by application of the 
imaging function $I: \mathbb{C}^N \to L^2(\Omega)$ defined by:
\begin{equation}
\label{def:I}
I[{\itbf a}](x,z) := 2 i \sum_{j=1}^N \beta_j e^{ i \beta_j x} {\phi_j(z)} \overline{a_j}  ,
\end{equation}
where  $\Omega\subset (a,+\infty) \times [0,L]$ is the compactly supported search region and the bar stands for complexe conjugate.
We can check that, in the case of a point-like source at $(x_o,z_o)$, $x_o > a$, 
we have 
$a_{j,o}= \frac{i}{2 \beta_j}  \phi_j(z_o) e^{i \beta_j x_o}$, $j=1,\ldots,N$, and
if we can estimate perfectly the vector ${\itbf a}_o$ from the data (which happens in particular when the antenna array spans the waveguide cross section, see below), then 
the imaging function has the form
\begin{equation}
\label{eq:expressaj0}
I[{\itbf a}_o](x,z) = \sum_{j=1}^N e^{ i \beta_j (x-x_o)} {\phi_j(z)} {\phi_j(z_o)} ,
\end{equation}
which is a peak centered at the source position $(x_o,z_o)$.
The resolution and stability properties of this imaging function (\ref{def:I}) have been analyzed in  \cite{bgt15}. 
The main result is that the width of the peak is approximately equal to the resolution limit $\lambda_o/2$, where
$\lambda_o=2\pi c_o/\omega$ is  the wavelength (with $c_o=$background velocity).

\begin{remark}
The imaging function (\ref{def:I}) is actually a reverse-time migration-type function \cite[Chapter 20]{book} (see also \cite{hodgkiss99,kuperman02,mordant,prada07,gar-pap-06}). 
Indeed, 
a reverse-time imaging function can be defined as $I_{\rm RT}:  L^2({\cal A},\mu) \to L^2(\Omega)$:
\begin{equation}
I_{\rm RT}[p](x,z) := -4  \int_{\cal A}   \partial_{x'}^2 \hat{G} (x,z;x',z')  \overline{p(x',z')} \mu(d{\itbf x}')   ,
\end{equation}
where $\hat{G}$ is the Green's function of the waveguide.
If we take into account only the guided modes, then the Green's function  has the form:
\begin{equation}
\hat{G} (x,z;x',z') = \frac{i}{2} \sum_{j=1}^N \frac{1}{\beta_j} e^{i \beta_j |x-x'|} \phi_j(z)  {\phi_j(z')}  ,
\end{equation}
and we find that, for $(x,z)\in \Omega$,
\begin{equation}
I_{\rm RT}[p](x,z)=  2i \sum_{j,j'=1}^N  {\beta_j} 
e^{ i \beta_j x }  {\phi_j(z)} \overline{ A_{jj'}} \overline{a_{j',o}}  ,
\end{equation}
where
\begin{equation}
A_{jj'}
=
\int_{{\cal A}}
e^{ i (\beta_{j}-\beta_{j'}) x' } {\phi_j(z')}  \phi_{j'}(z') \mu(d{\itbf x}')  ,
\label{def:prel}
\end{equation}
which is close to the  function $I[{\itbf a}_o](x,z)$ defined by (\ref{def:I})
when ${\bf A}$ is close to ${\bf I}$.
{Reverse-time migration functions are known to be efficient source imaging functions
as they can be seen as the solutions of least squares imaging \cite[Chapter 4, Section 4.1]{garnierpapa}.
They are the best estimators to localize point-like sources in the search domain (here, the interior of the waveguide),
in the sense that the position of the maximum of the modulus of the reverse-time migration function is the maximum likelihood estimator 
of the source position when the source is point-like and when the data are corrupted by additive noise~\cite{ammarigarniersolna}.}
\end{remark}

The imaging function (\ref{def:I}) is very efficient and has good resolution properties, 
but it requires to estimate the mode amplitudes ${\itbf a}_o$ of the 
recorded wavefield.
If the antenna array is dense, vertical and spans the full cross section of the waveguide, then 
the mode amplitudes ${\itbf a}_o$ can be easily obtained by projection of the observed wavefield $(p(x=0,z))_{z\in [0,L]}$ onto the mode profiles:
\begin{equation}
\label{def:aj}
\int_0^L p(x=0,z)  {\phi_j(z)} dz = a_{o,j} , \quad j=1,\ldots,N .
\end{equation}
We will see in the next section that it is possible to get good estimates of the mode amplitudes ${\itbf a}_o$ even when
the antenna array covers only a limited part of the cross section of the waveguide.

\section{Estimation of the mode amplitudes}
\label{sec:estim}
When the antenna array covers only a limited part  of the cross section of the waveguide
we would like to extract the vector ${\itbf a}_o$ from $(p(x,z))_{(x,z)\in {\cal A}}$ only.
This is actually possible, provided we know the mode profiles $(\phi_j(z))_{z \in [0,L]}$ and the modal wavenumbers $\beta_j$, $j=1,\ldots,N$.

\subsection{Perfect estimation}
In absence of any noise or measurement error, 
the following method can be implemented to estimate the vector ${\itbf a}_o$ 
(this is a general version of the weighted projection method proposed in \cite{tsogka16}):

\begin{enumerate}
\item 
Compute the Hermitian positive semi-definite matrix ${\bf A}$ of size $N\times N$ (as in (\ref{def:prel})):
\begin{equation}
\label{eq:defA}
A_{jl} := \int_{{\cal A}} \phi_j(z) \phi_l(z) e^{ i (\beta_{j}-\beta_{l}) x } \mu(d{\itbf x}) ,\quad j,l=1,\ldots,N .
\end{equation}
\item
Diagonalize the matrix 
$
{\bf A} = {\bf V}_{\! {\bf A}} {\bf D}_{\bf A} {\bf V}_{\! {\bf A}}^\dag  
$,
with ${\bf D}_{\bf A}$ diagonal matrix and ${\bf V}_{\! {\bf A}}$ unitary matrix
(here and below $\dag$ stands for conjugate and transpose).
\item
Introduce the reduced mode profiles:
\begin{equation}
\psi_l(x,z) := \sum_{j=1}^N  ({\bf V}_{\! {\bf A}})_{jl} \phi_j(z) e^{- i\beta_j x} ,\quad l=1,\ldots,N, \quad (x,z) \in {\cal A} .
\end{equation}
\item
 Compute the vector  $ {\itbf b} = ({b}_l)_{l=1}^N$ from the data $(p(x,z))_{(x,z) \in {\cal A}}$ by projection onto the reduced mode profiles:
\begin{equation}
\label{eq:projb}
b_l = \int_{{\cal A}} p(x,z) \overline{\psi_l}(x,z) \mu(d{\itbf x})  ,\quad l=1,\ldots,N  .
\end{equation}
\item
Compute the vector 
\begin{equation}
\label{eq:defab}
{\itbf a}    = {\bf V}_{\! {\bf A}} {\bf D}_{\bf A}^{-1} {\itbf b}.
\end{equation}
(If ${\bf A}$ is singular, then use the Moore-Penrose pseudo-inverse ${\bf D}_{\bf A}^+$ of ${\bf D}_{\bf A}$ instead of ${\bf D}_{\bf A}^{-1}$,
i.e. the diagonal matrix with 
diagonal coefficients $1/ ({\bf D}_{\bf A})_{jj}$ if $({\bf D}_{\bf A})_{jj}>0$ and $0$ otherwise).
\end{enumerate}

\begin{proposition}
If ${\bf A}$ is nonsingular, then ${\itbf a} ={\itbf a} _o$.
\end{proposition}

\noindent
{\it Proof.}
Let us study the method (\ref{eq:defA}-\ref{eq:defab}).
We have 
\begin{eqnarray*}
\int_{{\cal A}} \overline{\psi_l}(x,z) \phi_{l'}(z)  e^{-i\beta_{l'} x} \mu(d{\itbf x}) &=& \sum_{j=1}^N \overline{ ({\bf V}_{\! {\bf A}})_{jl}} \int_{{\cal A}} \phi_{j}(z) \phi_{l'}(z) e^{i(\beta_j-\beta_{l'})x } \mu(d{\itbf x})\\
& =&
 ( {\bf V}_{\! {\bf A}}^\dag {\bf A})_{ll'} 
= ( {\bf V}_{\! {\bf A}}^\dag {\bf V}_{\! {\bf A}} {\bf D}_{\bf A} {\bf V}_{\! {\bf A}}^\dag)_{ll'} = ( {\bf D}_{\bf A} {\bf V}_{\! {\bf A}}^\dag)_{ll'}   .
\end{eqnarray*}
 From  (\ref{eq:expressu0}), we get
$$
b_l = \sum_{l'=1}^N  \int_{{\cal A}} \overline{\psi_l}(x,z) \phi_{l'}(z) e^{-i\beta_{l'} x} \mu(d{\itbf x}) \, a_{l',o} 
= \sum_{l'=1}^N( {\bf D}_{\bf A} {\bf V}_{\! {\bf A}}^\dag)_{ll'}   a_{l',o}    ,
$$
i.e.,
$
{\itbf b} = {\bf D}_{\bf A}{\bf V}_{\! {\bf A}}^\dag {\itbf a}_o
$.
If ${\bf A}$ is nonsingular, then  ${\bf A}$ is positive definite and 
all eigenvalues of ${\bf D}_{\bf A}$ are not zero.
We then get  by (\ref{eq:defab}):
\begin{equation}
{\itbf a} =  {\bf V}_{\! {\bf A}} {\bf D}_{\bf A}^{-1} {\bf D}_{\bf A}{\bf V}_{\! {\bf A}}^\dag {\itbf a}_o= {\bf V}_{\! {\bf A}} {\bf V}_{\! {\bf A}}^\dag {\itbf a}_o ={\itbf a}_o ,
\end{equation}
the last equality follows from the unitarity of the matrix ${\bf V}_{\! {\bf A}}$.
\qed

\subsection{Regularized estimation}
The final step (\ref{eq:defab}) requires the matrix ${\bf A}$ to be positive-definite and well-conditioned for stability.
The conditioning of the matrix ${\bf A}$ is determined by the geometry of the array ${\cal A}$.
When the array does not cover the cross section of the waveguide, the conditioning of ${\bf A}$ may be poor and 
one should use a regularized pseudo-inverse for ${\bf D}_{\bf A}$~:
\begin{equation}
\label{def:aeps}
{\itbf a}_\epsilon  = {\bf V}_{\! {\bf A}} {\bf D}_{{\bf A}}^{\epsilon,+} {\itbf b}  ,
\end{equation}
where 
\begin{equation}
\label{def:D+eps}
{\bf D}_{{\bf A}}^{\epsilon,+} = {\rm Diag}\big( (\psi_\epsilon( ({\bf D}_{\bf A})_{jj}))_{j=1}^N \big)  ,
\end{equation}
 with 
\begin{equation}
\label{eq:psiepstyk}
 \psi_\epsilon(D_A) = D_A /(D_A^2 +\epsilon^2) \mbox{ (Tykhonov regularization)}
\end{equation}
 or 
\begin{equation}
{\psi}_\epsilon(D_A)=(1/D_A){\bf 1}_{(\epsilon,+\infty)}(D_A) \mbox{ (hard threshold regularization)}.
\end{equation}
We observe that we may not recover exactly the mode amplitudes when using the regularized method:
\begin{equation}
{\itbf a}_{\epsilon} =  
 {\bf V}_{\! {\bf A}} {\bf D}_{{\bf A}}^{\epsilon,+}  {\itbf b} =  {\bf V}_{\! {\bf A}}  {\bf D}_{{\bf A}}^{\epsilon,+} {\bf D}_{\bf A} {\bf V}_{\! {\bf A}}^\dag {\itbf a}_o 
 =
{\itbf a}_{o}
- 
{\bf V}_{\! {\bf A}} {\bf R}_{{\bf A}}^{\epsilon} {\bf V}_{\! {\bf A}}^\dag {\itbf a}_{o}  ,
\label{eq:biasreg}
\end{equation}
where the last term is an error term given in terms of  the diagonal matrix ${\bf R}_{{\bf A}}^{\epsilon}$  defined by
\begin{equation}
{\bf R}_{{\bf A}}^{\epsilon} = {\rm Diag}\big( (1-({\bf D}_{\bf A})_{jj} \psi_\epsilon ( ({\bf D}_{\bf A})_{jj}))_{j=1}^N \big).
\end{equation}
In the case of Tikhonov regularization, we have
$({\bf R}_{{\bf A}}^{\epsilon})_{jj} = \epsilon^2/(({\bf D}_{\bf A})_{jj}^2 +\epsilon^2) $.
In the case of hard threshold regularization, we have
$({\bf R}_{{\bf A}}^{\epsilon})_{jj} = {\bf 1}_{({\bf D}_{\bf A})_{jj}< \epsilon} $.

\subsection{Regularized estimation with measurement noise}
\label{subsec:regnoise}
As is well-known \cite{borcea10,tsogka16} and as is shown by (\ref{eq:biasreg}), regularization induces a bias, i.e. a deterministic error, but
it makes the estimation method much more robust with respect to noise, i.e. it can reduce the random error due to measurement noise.
This is a manifestation of the classical bias-variance tradeoff \cite{hastie}.
In order to illustrate this general statement, we here assume that the measurements 
$(p_{\rm meas}(x,z))_{z\in {\cal A}}$ are corrupted by an additive complex circular Gaussian noise:
\begin{equation}
p_{\rm meas}(x,z) 
=
p(x,z) + w(x,z) , \quad (x,z) \in {\cal A}  ,
\end{equation}
where $(w(x,z))_{(x,z)\in {\cal A}}$ is a Gaussian process with mean zero and 
delta covariance function:
\begin{equation}
\EE [w(x,z) \overline{w}(x',z') ] =\sigma^2
\left\{
\begin{array}{ll}
{\bf 1}_{z=z'} {\bf 1}_{x=x'} 
 ,&\mbox{ case (i),}\\
\displaystyle
\delta(z-z') \delta (x-x'),&\mbox{ case (ii),}\\
\displaystyle
\delta(z-z') ,&\mbox{ case (iii),}\\
\displaystyle
\delta (x-x'),&\mbox{ case (iv),} 
\end{array}
\right.
\end{equation}
for $(x,z),(x',z') \in {\cal A}$ (here ${\bf 1}_{z=z'}=1$ if $z=z'$ and $=0$ otherwise,
and $\delta$ is the Dirac distribution).

The estimated vector (\ref{def:aeps}) is here given by
\begin{equation}
{\itbf a}_\epsilon  = {\bf V}_{\! {\bf A}} {\bf D}_{{\bf A}}^{\epsilon,+} {\itbf b}_{\rm meas}  ,
\label{def:aepsbis}
\end{equation}
where the vector ${\itbf b}_{\rm meas} $ is obtained by projecting the measurements $(p_{\rm meas}(x,z))_{z\in {\cal A}}$
onto the reduced mode profiles as in (\ref{eq:projb}):
\begin{equation}
\label{def:bmeasl}
b_{{\rm meas},l} = \int_{{\cal A}} p_{\rm meas}(x,z) \overline{\psi_l}(x,z) \mu(d{\itbf x}) .
\end{equation}

\begin{proposition}
The mean square error consists of a bias term and a variance term:
\begin{eqnarray}
\EE \big[ \|{\itbf a}_\epsilon -{\itbf a}_o\|^2\big]
&=& 
\big\| \EE [ {\itbf a}_\epsilon  ] - {\itbf a}_o \big\|^2
+
\EE \big[ \|{\itbf a}_\epsilon -\EE [ {\itbf a}_\epsilon  ]\|^2\big] ,
\\
\big\| \EE [ {\itbf a}_\epsilon  ] - {\itbf a}_o \big\|^2
&=&
\sum_{j=1}^N [1-({\bf D}_{\bf A})_{jj}\psi_\epsilon(({\bf D}_{\bf A})_{jj})]^2 
|({\bf V}_{\! {\bf A}}^\dag {\itbf a}_{o})_{j}|^2  ,\\
\EE \big[ \|{\itbf a}_\epsilon -\EE [ {\itbf a}_\epsilon  ]\|^2\big] 
&=&
\sigma^2   \sum_{j=1}^N ({\bf D}_{\bf A})_{jj} \psi_\epsilon(({\bf D}_{\bf A})_{jj})^2 .
\end{eqnarray}
\end{proposition}
{\it Proof.}
The vector (\ref{def:bmeasl}) has the form
$$
b_{{\rm meas},l} = \int_{{\cal A}} p_{\rm meas}(x,z) \overline{\psi_l}(x,z) \mu(d{\itbf x}) 
  = \int_{{\cal A}} p(x,z) \overline{\psi_l}(x,z) \mu(d{\itbf x})  + w_l ,\quad l=1,\ldots,N  ,
$$
with $w_l = \int_{{\cal A}} w(x,z) \overline{\psi_l}(x,z) \mu(d{\itbf x})$. The random vector $(w_l)_{l=1}^N$ is Gaussian   with mean zero and covariance matrix: 
$$
\EE[ w_l \overline{w_{l'}}] = \sigma^2 \int_{{\cal A}} \overline{\psi_l}(x,z)  {\psi_{l'}}(x,z) \mu(d{\itbf x}) = \sigma^2 ({\bf V}_{\! {\bf A}}^\dag {\bf A} {\bf V}_{\! {\bf A}})_{ll'}=\sigma^2 ({\bf D}_{\bf A})_{ll'}  .
$$
This means that the random variables $w_l$ are independent Gaussian with mean zero and variances $\sigma^2 ({\bf D}_{\bf A})_{ll}$.

The estimated vector (\ref{def:aepsbis}) has mean 
$$
\EE [ {\itbf a}_\epsilon  ] ={\itbf a}_{o}
-  {\bf V}_{\! {\bf A}} {\bf R}_{\bf A}^\epsilon {\bf V}_{\! {\bf A}}^\dag {\itbf a}_{o}  ,
$$
and covariance
\begin{eqnarray*}
\EE \big[ ({\itbf a}_\epsilon -\EE [ {\itbf a}_\epsilon  ] ) 
( {\itbf a}_\epsilon - \EE [ {\itbf a}_\epsilon  ] )^\dag\big]
&=& 
\sigma^2 
{\bf V}_{\! {\bf A}}  {\bf D}_{{\bf A}}^{\epsilon,+} {\bf D}_{\bf A} {\bf D}_{{\bf A}}^{\epsilon,+}  {\bf V}_{\! {\bf A}}^\dag \\
&=&
\sigma^2 {\bf V}_{\! {\bf A}}
{\rm Diag}\big( (({\bf D}_{\bf A})_{jj}\psi_\epsilon(({\bf D}_{\bf A})_{jj})^2)_{j=1}^N \big) {\bf V}_{\! {\bf A}}^\dag  .
\end{eqnarray*}
The mean square error consists of a bias term and a variance term:
\begin{eqnarray*}
\EE \big[ \|{\itbf a}_\epsilon -{\itbf a}_o\|^2\big]
&=& 
\big\| \EE [ {\itbf a}_\epsilon  ] - {\itbf a}_o \big\|^2
+
\EE \big[ \|{\itbf a}_\epsilon -\EE [ {\itbf a}_\epsilon  ]\|^2\big]  ,
\end{eqnarray*}
with
\begin{eqnarray*}
\big\| \EE [ {\itbf a}_\epsilon  ] - {\itbf a}_o \big\|^2
&=&
\|{\bf V}_{\! {\bf A}} {\bf R}_{\bf A}^\epsilon {\bf V}_{\! {\bf A}}^\dag {\itbf a}_{o} \|^2  \\
&=&
\sum_{j=1}^N [1-({\bf D}_{\bf A})_{jj}\psi_\epsilon(({\bf D}_{\bf A})_{jj})]^2 
|({\bf V}_{\! {\bf A}}^\dag {\itbf a}_{o})_{j}|^2 ,
\end{eqnarray*}
and
\begin{eqnarray*}
\EE \big[ \|{\itbf a}_\epsilon -\EE [ {\itbf a}_\epsilon  ]\|^2\big]
&=&
\sigma^2 {\rm Tr}\big(  {\bf V}_{\! {\bf A}}
{\rm Diag}\big( (({\bf D}_{\bf A})_{jj}\psi_\epsilon(({\bf D}_{\bf A})_{jj})^2)_{j=1}^N \big) {\bf V}_{\! {\bf A}}^\dag \big) \\
&=&
\sigma^2   \sum_{j=1}^N ({\bf D}_{\bf A})_{jj} \psi_\epsilon(({\bf D}_{\bf A})_{jj})^2 .
\end{eqnarray*}
\qed

\begin{cor}
When $\sigma>0$, 
there exists a positive and finite $\epsilon$ that minimizes the mean square error.
\end{cor}
In other words,
regularization is always advantageous as soon as there is measurement noise.

\noindent
{\it Proof.}
For Tykhonov regularization (\ref{eq:psiepstyk}) the mean square error reads
$$
\EE \big[ \|{\itbf a}_\epsilon -{\itbf a}_o\|^2\big]
=
\sum_{j=1}^N \frac{\epsilon^4}{(({\bf D}_{\bf A})_{jj}^2 +\epsilon^2)^2}
|({\bf V}_{\! {\bf A}}^\dag {\itbf a}_{o})_{j}|^2
+
\sigma^2   \sum_{j=1}^N
\frac{ ({\bf D}_{\bf A})_{jj}^3}{(({\bf D}_{\bf A})_{jj}^2 +\epsilon^2)^2  }  .
$$
As $\epsilon \to 0^+$:
$$
\EE \big[ \|{\itbf a}_\epsilon -{\itbf a}_o\|^2\big]
=
\sigma^2   \sum_{j=1}^N
\frac{1}{({\bf D}_{\bf A})_{jj}}
-2\epsilon^2 \sigma^2   \sum_{j=1}^N 
\frac{1}{({\bf D}_{\bf A})_{jj}^3} +\dessous{O}{\epsilon \to 0}(\epsilon^4)  ,
$$
which shows that $\epsilon \in(0,+\infty) \mapsto \EE \big[ \|{\itbf a}_\epsilon -{\itbf a}_o\|^2\big]$
is a strictly decreasing function close to $0$.
As $\epsilon \to+\infty$:
$$
\EE \big[ \|{\itbf a}_\epsilon -{\itbf a}_o\|^2\big]
=
\sum_{j=1}^N 
|({\bf V}_{\! {\bf A}}^\dag {\itbf a}_{o})_{j}|^2
-
\epsilon^{-2}
\sum_{j=1}^N  ({\bf D}_{\bf A})_{jj}^2
|({\bf V}_{\! {\bf A}}^\dag {\itbf a}_{o})_{j}|^2
+
\dessous{O}{\epsilon \to +\infty}(\epsilon^{-4})  ,
$$
which shows that $\epsilon  \in(0,+\infty) \mapsto \EE \big[ \|{\itbf a}_\epsilon -{\itbf a}_o\|^2\big]$
is a strictly increasing function at infinity.
Since  $\epsilon  \in(0,+\infty) \mapsto \EE \big[ \|{\itbf a}_\epsilon -{\itbf a}_o\|^2\big]$ 
is continuous this shows that the exists an optimal $\epsilon \in (0,+\infty)$ that minimizes the mean square error and this optimal $\epsilon$ is positive and finite.
\qed

\section{Large dense antenna array}
\label{sec:large}
In this section we address the case of a large dense antenna array in a waveguide consisting of a large number of modes.
``Large antenna array" means much larger than the wavelength and ``dense antenna array" means that the Nyquist criterium is satisfied by the 
locations of the antennas so that the continuum approximation is valid.
``Large number of modes" means that the diameter of the cross section of the waveguide is much larger than the wavelength. 
This situation corresponds to a high-frequency regime (i.e. small wavelength),  which is not the main focus of this paper,
but this regime has motivated recent work in the literature.
We report in this section some interesting and original results about  the performances of
horizontal and vertical antenna arrays.

We will see below that the spectrum of the $N\times N$ matrix ${\bf A}$ corresponding to a large dense antenna array
typically contains two parts: $r_{\bf A}$ positive eigenvalues 
$({\bf D}_{\bf A})_{jj}$
for $j\leq r_{\bf A}$ and $N-r_{\bf A}$ vanishing eigenvalues 
$({\bf D}_{\bf A})_{jj} \simeq 0$ for $j>r_{\bf A}$. We can then say that $r_{\bf A}$ is the effective rank of the matrix,
and  the mean square error  is approximately:
$$
\EE \big[ \|{\itbf a}_\epsilon -{\itbf a}_o\|^2\big]
\simeq
\sum_{j=r_{\bf A}+1}^N  
|({\bf V}_{\! {\bf A}}^\dag {\itbf a}_{o})_{j}|^2  \simeq \frac{N-r_{\bf A}}{N} \|{\itbf a}_{o}\|^2,
$$
where we have used the rough approximation $|({\bf V}_{\! {\bf A}}^\dag {\itbf a}_{o})_{j}|^2  \simeq \| ({\bf V}_{\! {\bf A}}^\dag {\itbf a}_{o})\|^2 /N =\|{\itbf a}_{o} \|^2/N$. 
This shows that the quality of the estimation is directly related to the effective rank of the matrix ${\bf A}$ and the performance of the antenna array is all the better as its effective rank is larger.

\subsection{Vertical antenna array}
The case of a vertical antenna array occupying the line $\{0\}\times [0,2a]$ in a homogeneous waveguide
with bakground speed $c_o$ and Dirichlet boundary conditions is addressed in \cite{tsogka13,tsogka16}.
The number of guided modes is 
\begin{equation}
N = \lfloor \frac{k_o L }{ \pi} \rfloor = \lfloor  \frac{2L }{\lambda_o}\rfloor,
\end{equation}
 where $\lambda_o=2\pi c_o/\omega=2\pi/k_o$ is the wavelength.
In our framework, the problem is reduced to the analysis of the matrix
\begin{eqnarray}
\nonumber
A_{jl} &=& \frac{1}{2a} \int_0^{2a} \phi_j(z) \phi_l(z) dz \\
&=& 
\frac{1}{L} {\rm sinc}\big( \frac{2 \pi(l-j) a}{L}\big) -
\frac{1}{L} {\rm sinc}\big( \frac{2 \pi(l+j) a}{L}\big) ,
\label{eq:matAvert}
\end{eqnarray}
because $\phi_j(z) = \sqrt{2/L} \sin (\pi j z /L)$.
${\bf A}$ is a real, symmetric Toeplitz-minus-Hankel matrix.
Its spectral properties are determined by the Toeplitz part,
which can be studied in detail by the analysis conducted by Slepian about the discrete prolate spheroidal sequence \cite{slepian32}.
When $N \gg 1$ and $a/L =O(1)$, the spectrum can be decomposed into three parts:
there is a cluster of $O(N)$ eigenvalues close to $1/(2a)$, another cluster of $O(N)$  
eigenvalues close to $0$, and an intermediate layer of eigenvalues in between 
that decay from $1/(2a)$ to $0$. The number of eigenvalues in the intermediate layer is $o(N)$.
The number of ``significant" eigenvalues close to $1/(2a)$ is approximately equal to
$$
\Big[  N \frac{2a}{L}\Big] = \Big[ \frac{4a}{\lambda_o}\Big]  .
$$
The number of ``significant" eigenvalues, i.e. the effective rank of the matrix, is the length of the array $2a$ divided by the resolution limit $\lambda_o/2$.

The case of a vertical antenna array occupying the line $\{0\}\times [z_{\rm a}-a,z_{\rm a}+a]$ in a homogeneous waveguide
 is similar and the analysis of the previous case can be extended by the work of \cite{sengupta},
 as shown in \cite{tsogka16}. The results are similar in terms of numbers of ``significant" eigenvalues:
 The effective rank of the matrix ${\bf A}$ is the length of the array $2a$ divided by the resolution limit $\lambda_o/2$.

Finally, we can consider the general case of a vertical antenna array that occupies a set of disjoint intervals $\{0\}\times [b_k-a_k,b_k+a_k]$, $k=1,\ldots,P$
within $\{0\}\times (0,L)$. 

\begin{proposition}
The matrix ${\bf A}$ obtained with a vertical antenna array occupying the lines $\{0\}\times [b_k-a_k,b_k+a_k]$, $k=1,\ldots,P$, has an effective rank $2 N \sum_{k=1}^P a_k /L = 4 \sum_{k=1}^P a_k /\lambda_o$ when $N\to +\infty$.
\end{proposition}

\noindent
{\it Proof.}
The matrix ${\bf A}$ has the form
\begin{eqnarray*}
A_{jl} &=& \frac{1}{2 \sum_{k=1}^P a_k } \sum_{k=1}^P \int_{b_k-a_k}^{b_k+a_k} \phi_j(z) \phi_l(z) dz  \\
&=& \frac{1}{\pi} \int_0^\pi \big[ \cos((j-l)s)-\cos((j+l)s) \big] 
\rho(s) ds ,
\end{eqnarray*}
with
$$
\rho(s) =  \frac{1}{2 \sum_{k=1}^Pa_k }
\sum_{k=1}^P   {\bf 1}_{[\pi (b_k-a_k) / L, \pi (b_k+a_k)/ L]}(s)  .
$$
By \cite[Theorem 3.2]{fasino} the eigenvalues $(\sigma_j)_{j=1}^N$ of the matrix ${\bf A}$
satisfy for any continuous function $g$ 
$$
\frac{1}{N}
\sum_{j=1}^N g(\sigma_j)  \stackrel{N \to+\infty}{\longrightarrow} 
\frac{1}{\pi} \int_0^\pi g(\rho(s)) ds .
$$
This means that the empirical distribution of the eigenvalues $(\sigma_j)_{j=1}^N$ of the matrix ${\bf A}$
weakly converges as $N\to +\infty$ to a measure  supported by the two points $0$ and $[2\sum_{k=1}^P a_k]^{-1}$:
$$
\frac{1}{N}
\sum_{j=1}^N \delta_{\sigma_j} (d\sigma) \stackrel{N \to+\infty}{\longrightarrow} 
\Big(1- \frac{2 \sum_{k=1}^P a_k}{L} \Big) \delta_0(d\sigma) +    \frac{2 \sum_{k=1}^P a_k }{L} \delta_{1/[2\sum_{k=1}^P a_k ] }(d\sigma)  .
$$
This shows that the effective rank of the matrix is $2 N \sum_{k=1}^P a_k /L$.
\qed

In other words, 
the effective rank is the total length of the array $2 \sum_{k=1}^P a_k$ divided by the resolution limit $\lambda_o/2$.
It is interesting to note that, for a homogeneous waveguide and in the continuum approximation, 
the spatial distribution of the receivers along the vertical cross section does not play any role, only the 
total length of the linear antenna array plays a role.
 
 \begin{figure}
\begin{center}
\begin{tabular}{c}
{\bf a)}
\includegraphics[width=6.cm]{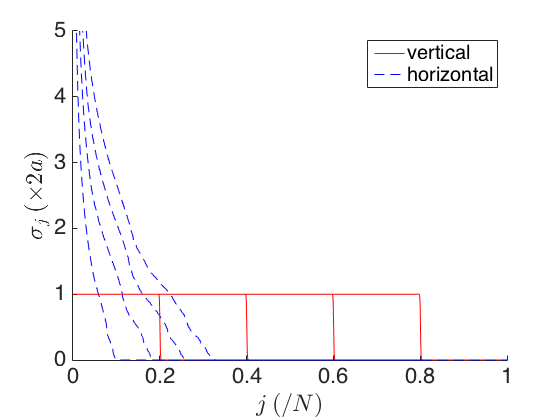} 
{\bf b)}
\includegraphics[width=6.cm]{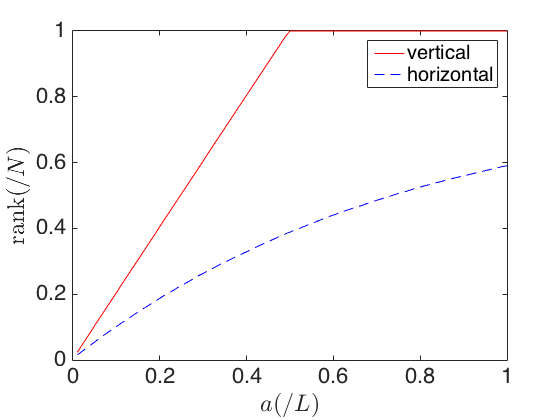} 
\end{tabular}
\end{center}
\caption{Picture a: Eigenvalues $(\sigma_j)_{j=1}^N$ of the matrix ${\bf A}$ for different values of the ratio $a/L$ ($0.1$, $0.2$, $0.3$, $0.4$)
and for a vertical array (red solid) and for a horizontal array (blue dashed).
Here $L=1000$, $k_o=1$, $z_{\rm a}=220$.
The vertical array has about $N \times 2a/L$ significant eigenvalues with value $1/(2a)$.
The horizontal array with the same length has less significant eigenvalues.
Picture b: Effective rank of the matrix ${\bf A}$ for different values of the ratio $a/L$. The matrix is full rank for a vertical array with length $2a=L$ which covers the whole cross section.}
\label{fig:verthori}
\end{figure}

\subsection{Horizontal antenna array}
The case of a horizontal antenna array occupying the line $[0,2a] \times \{z_{\rm a}\}$ in a homogeneous waveguide
 is qualitatively similar. The matrix
${\bf A}$ has the form
\begin{eqnarray}
\nonumber
A_{jl} &=& \frac{1}{2a} \int_{0}^{2a} \phi_j(z_{\rm a}) \phi_l(z_{\rm a}) e^{i (\beta_j-\beta_l)x} dx \\
&=& 
\frac{2}{L}\sin \big(\frac{\pi j z_{\rm a}}{L}\big)e^{i \beta_j a}  \sin\big( \frac{\pi l z_{\rm a}}{L}\big) e^{-i \beta_l a}  {\rm sinc} \big((\beta_j-\beta_l) a\big) .
\label{eq:matAhor}
\end{eqnarray}
Unless $z_{\rm a}$ corresponds to a node of a mode, the spectral properties of ${\bf A}$ are related to those
of the matrix $\tilde{\bf A}=\big( {\rm sinc} \big((\beta_j-\beta_l) a\big)\big)_{j,l=1}^N$,
which looks like the sinc kernel addressed by Slepian, upon substitution $\pi j /L \mapsto \beta_j=\sqrt{k_o^2 - \pi^2 j^2/L^2}$.
The theoretical analysis of this case, as far as we know, has not yet been carried out.
We will first do numerical simulations to propose some conjectures 
and then we will give the theoretical results.

Based on numerical simulations (see figure \ref{fig:verthori}), we get the following conjecture:
When $N \gg 1$ and $a/L =O(1)$, the spectrum can be decomposed into two parts:
there is a cluster of eigenvalues of order $1/a$ and
another cluster of eigenvalues close to $0$  (see figure \ref{fig:verthori}a).
The number of significant eigenvalues is approximately equal to 
$ [  N {a}/{L}  ]$ when $a/L$ is small, and smaller than $ [  N {a}/{L}  ]$ when $a/L$ becomes of order one (see figure \ref{fig:verthori}b).
Note that $ [  N {a}/{L}  ]$ is one half the number of significant eigenvalues 
for a vertical antenna array with the same length. This conjecture is proved in the following proposition in a more general case.

We now address the case of a horizontal antenna array occupying the disjoint intervals $[b_k-a_k,b_k+a_k] \times \{z_{\rm a}\}$, $k=1,\ldots,P$, in a homogeneous waveguide,
with $\cup_{k=1}^P [b_k-a_k,b_k+a_k] \subset [-L,L]$.

\begin{proposition}
For almost every $z_{\rm a} \in (0,L)$, 
the matrix ${\bf A}$ obtained with a horizontal antenna array occupying the lines $[b_k-a_k,b_k+a_k] \times \{z_{\rm a}\}$, $k=1,\ldots,P$, has an effective rank equal to $N \sum_{k=1}^P a_k /L =2\sum_{k=1}^P a_k /\lambda_o $
when $N \to +\infty$ and the total length of the antenna array is much smaller than $L$.
\end{proposition}

\noindent
{\it Proof.}
The matrix
${\bf A}$ has the form
\begin{eqnarray}
\nonumber
A_{jl} 
&=& \frac{1}{2\sum_{k=1}^P a_k}
\sum_{k=1}^P  \int_{b_k-a_k}^{b_k+a_k} \phi_j(z_{\rm a}) \phi_l(z_{\rm a}) e^{i (\beta_j-\beta_l)x} dx \\
&=& 
\frac{1}{L}\sin \big(\frac{\pi j z_{\rm a}}{L}\big) \sin\big( \frac{\pi l z_{\rm a}}{L}\big)  
 \frac{1}{\sum_{k=1}^P a_k} 
\sum_{k=1}^P 2a_k e^{i (\beta_j-\beta_l) b_k} {\rm sinc} \big((\beta_j-\beta_l) a_k\big) .
\end{eqnarray}

We first show the following result: If ${\itbf U}$ is a vector with non-zero entries
and $A_{jl}=U_j \tilde{A}_{jl} U_l$ is a $N\times N$ symmetric real matrix, then the rank of ${\bf A}$ and $\tilde{\bf A}$ are equal. 
Indeed, if $\tilde{r}$ is the rank of $\tilde{\bf A}$, then $
\tilde{\bf A} = \sum_{k=1}^{\tilde{r}} \tilde{\sigma}_k \tilde{\itbf v}_k\tilde{\itbf v}_k^T$
with orthornomal vectors $\tilde{\itbf v}_k$ and nonzero $\tilde{\sigma}_k$,
and therefore
${\bf A} = \sum_{k=1}^{\tilde{r}} \tilde{\sigma}_k  {\itbf v}_k {\itbf v}_k^T$ with ${v}_{k,j}= \tilde{v}_{k,j} U_j$
for $j=1,\ldots,N$ and $k=1,\ldots,\tilde{r}$.
The vectors $ {\itbf v}_k $ are linearly independent 
(if $\sum_{k=1}^{\tilde{r}} \alpha_k {\itbf v}_k={\bf 0}$, then $\sum_{k=1}^{\tilde{r}} \alpha_k \tilde{\itbf v}_k={\bf 0}$, 
and therefore $\alpha_k=0$ for all $k$).
This shows that the rank of ${\bf A}$ is $\tilde{r}$.

From the previous result, for almost every $z_{\rm a} \in (0,L)$ (for all $z_{\rm a}$ except $z_{\rm a} \in \{ L/k, \, k\in \{2,3,4,\ldots\}\}$, so that $\sin (\pi j z_{\rm a}/L)$ never cancels), the matrix ${\bf A}$ has the same rank as the matrix $\tilde{\bf A} $ with
\begin{eqnarray*}
\tilde{A}_{jl} = {\cal K}\big(\frac{(\beta_j-\beta_l)L}{\pi}\big),\quad \quad 
{\cal K}(s) =  \frac{1}{2\sum_{k=1}^P a_k} 
\sum_{k=1}^P 2a_k e^{i \pi  s \frac{b_k}{L}} {\rm sinc} \big( \pi  s \frac{a_k}{L}\big) .
\end{eqnarray*}
We have
\begin{eqnarray*}
(\tilde{\bf A} {\itbf u})_j &=& N \int_0^1 {\cal K}\Big(\frac{(\beta_j-\beta_{\lceil sN \rceil})L }{\pi }\Big)   u_{\lceil sN \rceil} ds \\
&=& N \int_0^1 {\cal K}\Big(\frac{\beta_j L}{\pi } - N \sqrt{1-\frac{\lceil sN\rceil^2}{N^2}} \Big)   u_{\lceil sN \rceil} ds  .
\end{eqnarray*}
When the length of the antenna array is much smaller than $L$, then we can make the continuous approximation ${\cal K}\big(\frac{\beta_j L}{\pi} - N \sqrt{1-\frac{\lceil s N\rceil^2}{N^2}} \big) \simeq {\cal K}\big(\frac{\beta_j L}{\pi} - N \sqrt{1-s^2} \big)$.
Therefore
\begin{eqnarray*}
(\tilde{\bf A} {\itbf u})_j 
&=& N \int_0^1 {\cal K}\Big(\frac{\beta_j L}{\pi } - N \sqrt{1-s^2} \Big)   u_{\lceil s N\rceil} ds 
\\
&=& N \int_0^1 {\cal K}\Big(\frac{\beta_j L}{\pi } - N s'  \Big)   u_{\lceil\sqrt{1-s'^2}N\rceil} \frac{s'}{\sqrt{1-s'^2}} ds'  ,
\end{eqnarray*}
and
\begin{eqnarray*}
(\tilde{\bf A} {\itbf u})_{\lceil\sqrt{1-s^2}N\rceil} = N \int_0^1 {\cal K}\big(N(s- s')  \big)   u_{\lceil\sqrt{1-s'^2}N\rceil} \frac{s'}{\sqrt{1-s'^2}} ds' .
\end{eqnarray*}
If we introduce $\tilde{\itbf u}_j= \frac{\sqrt{j}}{\sqrt[4]{1+N^2-j^2}} {\itbf u}_{\lfloor \sqrt{1+N^2-j^2} \rfloor}$, we get
\begin{eqnarray*}
(\tilde{\bf A} \tilde{\itbf u})_{\lceil s N \rceil} = N \int_0^1 \frac{\sqrt{s}}{\sqrt[4]{1-s^2}}  {\cal K}\big(N(s- s')  \big) \frac{\sqrt{s'}}{\sqrt[4]{1-s'^2}}  \tilde{u}_{\lceil s' N\rceil}  ds'  .
\end{eqnarray*}
Using the same argument as above (multiplication left and right by the same vector does not change the rank), 
we conclude that the rank of $\tilde{\bf A}$ is the same as the rank of
the Toeplitz matrix $({\cal T}( j-l))_{j,l=1}^N$:
$$
{\cal T}(j-l) =  \frac{1}{2\pi} \int_{-\pi}^\pi  e^{ i (j-l)s} 
\rho(s) ds ,
$$
with
$$
\rho(s) =  \frac{1}{2 \sum_{k=1}^Pa_k }
\sum_{k=1}^P   {\bf 1}_{[\pi (b_k-a_k) / L, \pi (b_k+a_k)/ L]}(s)  .
$$
The rank of $({\cal T}( j-l))_{j,l=1}^N$  is $N \sum_{k=1}^P a_k /L $ when $N \to +\infty$ by \cite[Section 5.2]{grenander}.
\qed

To summarize, the results for the large dense antenna array show that the vertical arrays perform better than horizontal arrays (with the
same lengths). The length of the horizontal array should be twice as long as the one of the vertical array to present similar performance
(in the sense that the same amount of information can be extracted from the two arrays).

\section{Small discrete antenna array}
\label{sec:small}
We consider in this section that
 the antenna array is discrete and consists of  $M$ point-like receivers  localized at $(x_k,z_k)$, $k=1,\ldots,M$.
Then the recorded signals are (for sources located in the region $x>a$):
\begin{equation}
p_k = p(x_k,z_k) = \sum_{j=1}^N a_{j,o} \phi_j(z_k) \exp (- i \beta_j x_k)  ,\quad
k=1,\ldots,M.
\end{equation}
The recorded vector ${\itbf p}=(p_k)_{k=1}^M$
has the form
\begin{equation}
{\itbf p} = {\bf B} {\itbf a}_o  ,
\end{equation}
where ${\bf B}$ is the $M \times N$ matrix with entries 
\begin{equation}
\label{def:matB}
B_{kj}=  \phi_j(z_k)  \exp (-i \beta_j x_k), \quad k=1,\ldots,M, \quad j=1,\ldots,N.
\end{equation}
Throughout the section we assume that $M \geq N$ (i.e. there are more receivers than guided modes).
The matrix ${\bf A}$ in (\ref{eq:defA}) has the form
\begin{equation}
{\bf A} = \frac{1}{M} {\bf B}^\dag {\bf B},
\end{equation}
which shows that the singular values of ${\bf B}$ are the square roots of the eigenvalues of ${\bf A}$ (up to the factor $1/M$)
and the right singular vectors of ${\bf B}$ are the eigenvectors of ${\bf A}$ (i.e. the columns of the matrix ${\bf V}_{\! {\bf A}}$).
Therefore it is possible to work directly with the inverse problem associated with~${\bf B}$ in the case of discrete antenna arrays.

\subsection{Source imaging function}
As explained in Section \ref{sec:imag}, source imaging has two main steps:\\
1) regularized estimation of the mode amplitudes,\\
2) migration of the estimated mode amplitudes by the imaging function (\ref{def:I}).

It is possible to recover all mode amplitudes ${\itbf a}_o$ from the vector of recorded signals ${\itbf p}$, provided 
${\bf B}$ has rank $N$. 
The ideal method consists in applying the pseudo-inverse of ${\bf B}$ to the observed vector ${\itbf p}$:
\begin{equation}
{\itbf a} = {\bf B}^+ {\itbf p}  ,
\end{equation}
where ${\bf B}^+={\bf V} {\bf D}^+ {\bf U}^\dag$, ${\bf D}^+$ is the $N\times M$ diagonal matrix with 
diagonal coefficients $1/D_{jj}$ if $D_{jj}>0$ and $0$ otherwise, and 
\begin{equation}
{\bf B} = {\bf U}  {\bf D} {\bf V}^\dag
\end{equation}
 is the singular value decomposition of ${\bf B}$.
We then have
${\itbf a} = {\bf V} {\bf D}^+{\bf D} {\bf V}^\dag {\itbf a}_o$,
which shows that, if ${\bf B}$ has rank $N$, then ${\bf D}^+{\bf D}={\bf I}$ and 
${\itbf a} ={\itbf a}_o$.

In practice, one needs to use a regularized pseudo-inverse ${\bf D}^+_\epsilon$ instead of ${\bf D}^+$ as in (\ref{def:D+eps}).
This has to be done in particular when there is measurement noise.
The regularization discards the contributions that correspond to small singular values,
because they cannot be estimated with accuracy.
If the recorded vector ${\itbf p}_{\rm meas}$ has the form
\begin{equation}
{\itbf p}_{\rm meas} = {\itbf p}  +{\itbf w} , \quad \quad {\itbf p}= {\bf B} {\itbf a}_o,
\end{equation}
with 
\begin{equation}
 {\itbf w}=(w_j)_{j=1}^M \sim {\cal N}({\bf 0}, \sigma_{\rm meas}^2{\bf I}),
\end{equation}
i.e.,  
a family of independent and identically distributed Gaussian circular complex
random variables with mean zero and variance $\sigma_{\rm meas}^2$,
then the vector recovered by application of the regularized pseudo-inverse ${\bf B}^+_\epsilon=
 {\bf V} {\bf D}^+_\epsilon  {\bf U}^\dag$ 
is
\begin{equation}
{\itbf a}_\epsilon =  {\bf B}^+_\epsilon {\itbf p}_{\rm meas}
=
 {\bf V} {\bf D}^+_\epsilon {\bf D} {\bf V}^\dag {\itbf a}_o + {\bf V} {\bf D}^+_\epsilon {\bf U}^\dag {\itbf w} .
\end{equation}
Since ${\bf U}^\dag {\itbf w} \sim {\cal N}({\bf 0}, \sigma_{\rm meas}^2{\bf I})$,
 the mean square error of the 
Tykhonov-regularized estimator  is the sum of a bias term and a variance term:
\begin{equation}
\EE \big[ \| {\itbf a}_\epsilon  - {\itbf a}_o\|^2\big]
=
\sum_{j=1}^N \frac{\epsilon^4}{(D_{jj}^2 +\epsilon^2)^2 }  \big|({\bf V}^\dag {\itbf a}_o )_j \big|^2
+
\sum_{j=1}^N \frac{D_{jj}^2}{(D_{jj}^2 +\epsilon^2)^2 }  \sigma_{\rm meas}^2 .
\end{equation}
To be complete, we can remark that this regularization corresponds to the choice $\psi^\epsilon(D_A) = 1/(D_A+\epsilon^2)$
in the general framework of Section \ref{subsec:regnoise}, because $D_{jj}=({\bf D}_{\bf A})_{jj}^{1/2}$.

As in Section \ref{subsec:regnoise} one can show that it is always advantageous to regularize.
Indeed the mean square error is a smooth function of the regularization parameter $\epsilon$, 
it is strictly decreasing close to zero and strictly increasing at infinity.
The optimal regularization parameter satisfies $\partial_{(\epsilon^2)} \EE \big[ \| {\itbf a}_\epsilon  - {\itbf a}_o\|^2\big]=0$, that is to say,
\begin{equation}
\sum_{j=1}^N \frac{D_{jj}^2}{(D_{jj}^2 +\epsilon^2)^3 }\big(  \big| ({\bf V}^\dag {\itbf a}_o )_j \big|^2 \epsilon^2 -   \sigma_{\rm meas}^2\big) =0 .
\end{equation}
By using the rough approximation $ \big| ({\bf V}^\dag {\itbf a}_o )_j \big|^2 
\simeq \|  {\bf V}^\dag {\itbf a}_o  \|^2/N=\|  {\itbf a}_o \|^2/N$, 
this shows that 
\begin{equation}
\epsilon^2 \simeq \sigma_{\rm meas}^2  \frac{N}{\|{\itbf a}_o\|^2} =  \frac{\sigma_{\rm meas}^2}{\frac{1}{N} \sum_{j=1}^N |a_{o,j}|^2} .
\end{equation}
In other words, the regularization parameter $\epsilon$ should be proportional to the standard deviation of the measurement error.
This is a standard choice for the Tikhonov regularization parameter \cite{scherzer}.
This choice is also promoted by the Morozov's discrepancy principle, 
which claims that we should not try to fit the data beyond the measurement noise \cite{engl}.

The quality of the image {$(x,z)\mapsto |I[{\itbf a}_\epsilon](x,z)|$}  built from the estimation of the vector~${\itbf a}_\epsilon$ using the imaging 
function (\ref{def:I}) depends on the noise level and the conditioning
of the matrix ${\bf B}$, which itself depends on the array geometry.
In the following subsections we analyze different array geometries and determine the effective rank
of the matrix ${\bf B}$.

\subsection{Vertical antenna array}
 When the vertical antenna array in the plane $x=0$ consists of  $M$ point-like receivers located at $z_k$, 
$k=1,\ldots,M$, the matrix ${\bf B}$ has the form
\begin{equation}
{\bf B} =  \big( \phi_j(z_k) \big)_{1\leq k \leq M , 1\leq j \leq N} .
\end{equation}
Let us consider the case where the velocity is constant and equal to $c_o$ and the two boundary conditions are Dirichlet.
Then 
\begin{equation}
\label{eq:expressmodes}
\phi_j(z) = \frac{\sqrt{2}}{\sqrt{L}} \sin \big( \alpha_j z\big), \quad 
\quad \alpha_j =  \frac{\pi j  }{L}.
\end{equation}
By denoting by $a$ the radius of the antenna array, by $z_{\rm a}$ its center, and by introducing $z_k=z_{\rm a} +\tilde{z}_k$
(so that $|\tilde{z}_k|\leq a$),
the matrix ${\bf B}$ can be expanded as
\begin{equation}
{\bf B} =
\sum_{q=0}^{Q-1}
{\itbf u}_q {\itbf v}_q^T
+O\Big( \frac{ (k_oa)^{Q}}{Q!} \Big)   ,
\label{eq:expandB1}
\end{equation}
when $k_o a\ll 1$ ($k_o=\omega/c_o=2\pi/\lambda_o$ is the homogeneous wavenumber), 
with
\begin{eqnarray}
{\itbf u}_{2q+1} &=& \Big( \frac{\sqrt{2}}{\sqrt{L}} \cos(\alpha_j z_{\rm a})  \alpha_j^{2q+1} \Big)_{j=1}^N,\\
{\itbf u}_{2q} &=&   \Big( \frac{\sqrt{2}}{\sqrt{L}}\sin(\alpha_j z_{\rm a}) \alpha_j^{2q} \Big)_{j=1}^N,\\
{\itbf v}_{2q+1} &=& \frac{(-1)^{q}}{(2q+1)!}
\Big( \tilde{z}_k^{2q+1}\Big)_{k=1}^M,\\
{\itbf v}_{2q} &=& \frac{(-1)^q}{(2q)!}
\Big( \tilde{z}_k^{2q}\Big)_{k=1}^M.
\end{eqnarray}

\begin{proposition}
\label{prop:linind1}
The vectors $({\itbf v}_q)_{q=0}^{Q-1}$ are linearly independent when $Q\leq M$.
The vectors $({\itbf u}_q)_{q=0}^{Q-1}$ are linearly independent when $Q\leq N$,
except possibly for a finite number of special values of $z_{\rm a}$.
\end{proposition}

\noindent
{\it Proof.}
If $\sum_{q=0}^{Q-1} \lambda_q {\itbf v}_q={\bf 0}$, then the polynomial $z \mapsto \sum_{q=0}^{Q-1} \lambda_q \frac{(-1)^{[q/2]}}{q!}
z^q$ has $M$ distinct zeros $(\tilde{z}_k)_{k=1}^M$. If $M>Q-1$, then the polynomial must be zero which imposes $\lambda_q=0$ for all $q$.\\
If $\sum_{q=0}^{N-1} \lambda_q {\itbf u}_q={\bf 0}$ has a non trivial solution $(\lambda_q)_{q=0}^{N-1}$, then the determinant
of the matrix $ (\partial_{z_{\rm a}}^{q-1} \sin (\alpha _j z_{\rm a}) )_{q,j=1}^N$ is zero. By Euler's formula, the function 
$z \mapsto {\rm det} \big((\partial_z^{q-1} \sin (\alpha _j z))_{q,j=1}^N\big)$ can 
be written as $\exp(- i N(N+1) \pi z/(2L)) Q_{N}(\exp ( i  \pi z /L))$ where $Q_{N}$ is a polynomial of degree $N(N+1)$.
This non-zero polynomial has only a finite number of roots, 
so there is only a finite number of $z\in [0,L]$ such that $Q_{N}(\exp ( i N \pi z /L)) =0$.
If $z_{\rm a}$ is different from these special values, then the vectors $({\itbf u}_q)_{q=0}^{N-1}$ are linearly independent.
\qed

The $\epsilon$-regularization  used when the noise level has relative standard deviation $\epsilon$
prevents from exploiting the singular vectors whose singular values are smaller than 
$\epsilon$. 
From Eq.~(\ref{eq:expandB1}) and Proposition \ref{prop:linind1} this implies that ${\bf B}$ has an effective rank $Q$
where $Q$ is such that 
\begin{equation}
\label{eq:condeps}
\frac{ (k_oa)^{Q}}{Q!}\simeq \epsilon .
\end{equation}
Here we assume that $M$ is larger than $Q$ (otherwise the rank is limited to $M$).
For instance, for $\epsilon=10^{-7}$ and $k_oa=0.125$, formula (\ref{eq:condeps}) predicts that 
we have $Q\simeq 5$.
If we compare with the singular values of the matrix ${\bf B}$ when $k_o =1 $, $L=20$, $z_k=11+0.25 (k-M/2)/M$, $k=1,\ldots,M$, $M=20$,
then we find that $N=6$ and the first singular values of ${\bf B}$ are $\sigma(1)\simeq2.6$, $\sigma(2)\simeq 0.1$,
$\sigma(3)=2\, 10^{-3}$, $\sigma(4)=1 \, 10^{-5}$, $\sigma(5)=1 \, 10^{-7}$,
$\sigma(6)=4 \, 10^{-10}$,
so that indeed its effective rank is approximately $5$.

\begin{figure}
\begin{center}
\begin{tabular}{cc}
\includegraphics[width=6.cm]{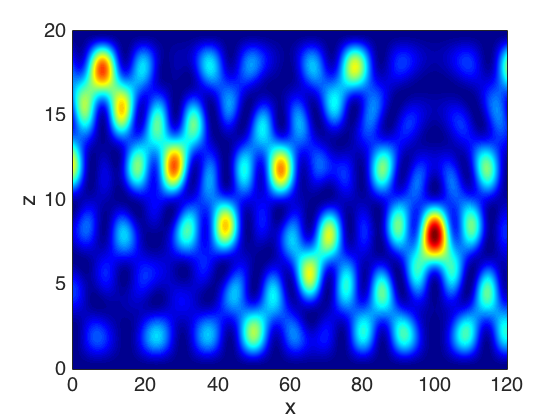}&
\includegraphics[width=6.cm]{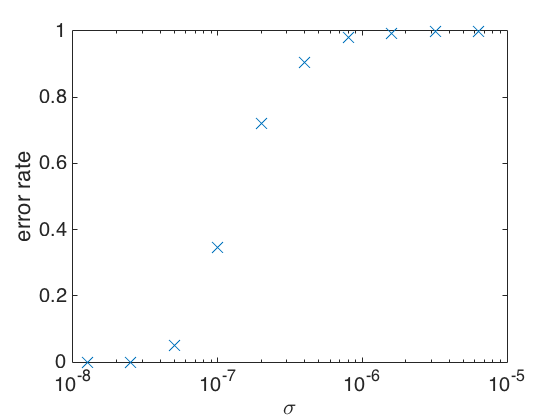}
\\
$\sigma=0$ & 
localization error rate
 \\
\includegraphics[width=6.cm]{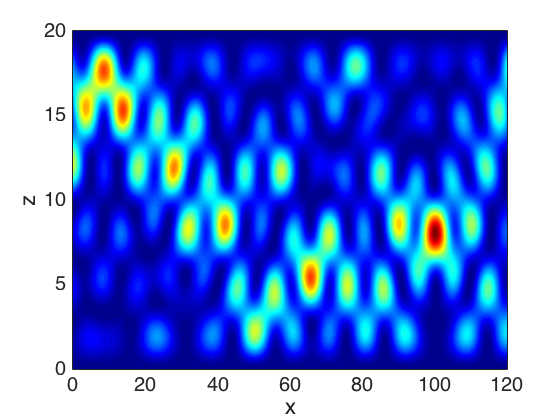}&
\includegraphics[width=6.cm]{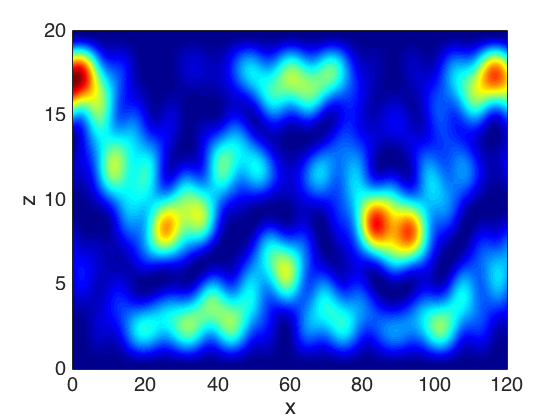}\\
$\sigma=10^{-7}$ & $\sigma=10^{-6}$
\end{tabular}
\end{center}
\caption{
{
Images and localization error rates obtained with a vertical antenna array of $M=20$ receivers 
and total length $0.25 \simeq 0.04 \lambda_o$ and with different levels of noise~$\sigma$. Here $\sigma_{\rm meas} = \sigma \| {\itbf p}\|_\infty$.
The source position is at $(x_o,z_o)=(100,7.7)$. The source can be localized with accuracy (at the scale of the wavelength) 
if $\sigma \lesssim 10^{-7}$.
}
}
\label{fig:1}
\end{figure}

In figure \ref{fig:1},
a vertical antenna array records the time-harmonic wave. It is localized at $z_{\rm a}=11$ and $z_k=11+0.25 (k-M/2)/M$, $k=1,\ldots,M$,
 with $M=20$.
Here the frequency is $\omega=1$, the velocity is $c_o=1$, $L=20$, the original source is at $(x_o,z_o)=(100,7.7)$.
There are $N=6$ guided modes.
Different noise levels are considered (corresponding to different regularization parameters). 
{The imaging function $(x,z) \mapsto |I[{\itbf a}_\epsilon](x,z)|$ is plotted for $(x,z)$ within the waveguide for different values of the noise level. The imaging function is normalized by its maximal value. We can observe that the position of the maximum of the imaging function corresponds to 
the source position with very high probability and with very good accuracy when the noise level is small. 
There exists a critical noise level beyond which the method fails, the imaging function has many local maxima and
 the position of the global maximum of the imaging function does not correspond anymore to the source position.  
 We also plot in  figure \ref{fig:1} the localization error rate as a function of the noise level, it is the probability that the position of the maximum of the imaging function is  less than half-a-wavelength away from the source position, it is computed by an empirical average based on $1000$ simulations with independent and identically distributed noise realizations.
 The numerical results show that the source can be localized with accuracy (at the scale of the wavelength)  and with high probability
if $\sigma \lesssim 10^{-7}$.}
Note that the total length of the array is very small compared to the wavelength, it is equal to
$0.25\simeq 0.04 \lambda_o$. This shows that it is possible to image the source with such an antenna array, but the signal-to-noise
ratio has to be very high.

\subsection{Horizontal antenna array}
In this subsection we consider the situation in which the antenna array is horizontal at $z=z_{\rm a}$ 
and consists of $M$ receivers localized at $x=x_k$, $k=1,\ldots,M$, around the position $x=0$.
The matrix ${\bf B}$ has the form
\begin{equation}
{\bf B}= \big(  \phi_j(z_{\rm a})  \exp (-i \beta_j x_k) \big)_{1\leq k\leq M, 1\leq j \leq N}.
\end{equation}

Let us consider the case where the velocity is constant and the two boundary conditions are Dirichlet.
Then the mode profiles are given by (\ref{eq:expressmodes}) and the modal wavenumbers are
\begin{equation}
\label{eq:expressmodes2}
\beta_j  =  \sqrt{k_o^2 - \alpha_j^2} .
\end{equation}
By denoting by $a$ the length of the antenna array, 
${\bf B} $ can be expanded as
\begin{equation}
{\bf B}= 
\sum_{q=0}^{Q-1}
{\itbf u}_q {\itbf v}_q^\dag
+O\Big( \frac{ (k_oa)^{Q}}{Q!} \Big)
\end{equation}
when $k_o a\ll 1$, 
with
\begin{equation}
{\itbf u}_q = \Big( \phi_j(z_{\rm a}) \beta_j^{q}\Big)_{j=1}^N,
\quad 
{\itbf v}_q = \frac{i^q}{q!} \Big( x_k^{q}\Big)_{k=1}^M.
\end{equation}

\begin{proposition}
The vectors $({\itbf v}_q)_{q=0}^{Q-1}$ are linearly independent when $Q\leq M$.
The vectors $({\itbf u}_q)_{q=0}^{Q-1}$ are linearly independent when $Q\leq N$,
except possibly for a finite number of special values of $z_{\rm a}$.
\end{proposition}
{\it Proof.}
If $\sum_{q=0}^{Q-1} \lambda_q {\itbf v}_q={\bf 0}$, then the polynomial $z \mapsto \sum_{q=0}^{Q-1} \lambda_q \frac{i^q}{q!}
x^q$ has $M$ distinct zeros $(x_k)_{k=1}^M$. If $M>Q-1$, then the polynomial must be zero which imposes $\lambda_q=0$ for all~$q$.\\
If $\sum_{q=0}^{N-1} \lambda_q {\itbf u}_q={\bf 0}$ has a non trivial solution $(\lambda_q)_{q=0}^{N-1}$, then the determinant
of the matrix $ ( \beta_j^{q-1} \sin (\alpha_j z_{\rm a}) )_{q,j=1}^N$ is zero. The function 
$z \mapsto {\rm det} \big((\beta_j^{q-1}  \sin (\alpha_j z))_{q,j=1}^N\big)$ can 
be written as $\exp(- i N(N+1) \pi z/(2L)) Q_{N}(\exp ( i \pi z /L))$ where $Q_{N}$ is a polynomial of degree $N(N+1)$.
Therefore there is only a finite number of $z\in [0,L]$ such that $Q_{N}(\exp ( i \pi z /L)) =0$.
If $z_{\rm a}$ is different from these special values, then the vectors $({\itbf u}_q)_{q=0}^{N-1}$ are linearly independent.
\qed

The $\epsilon$-regularization used when the noise level has relative standard deviation $\epsilon$
prevents from exploiting the singular vectors whose singular values are smaller than 
$\epsilon$. 
This implies that ${\bf B} $ has an effective rank $Q$
where $Q$ is such that $\frac{ (k_oa)^{Q}}{Q!}\simeq \epsilon$.
Here we assume that $M$ is larger than $Q$.
For instance, if $\epsilon =10^{-7}$ and $k_oa=0.125$, then $Q\simeq 5$.
If we compare with the singular values of the matrix ${\bf B}$ when $k_o =1 $, $L=20$, $z_{\rm a}=11$, $x_k=0.25(k-M/2)/M$, $k=1,\ldots,M$, $M=20$,
then we find that the first singular values of ${\bf B}$ are $\sigma(1)\simeq 2.6$, $\sigma(2)\simeq 0.04$,
$\sigma(3)=3 \, 10^{-4}$, $\sigma(4)=7 \, 10^{-7}$, $\sigma(5)=8 \, 10^{-10}$,
so that indeed its effective rank is $4$. We observe that the singular values decay slightly faster than in the case of a vertical array.
This is due to the fact that the $\beta_j$ are not uniformly distributed over $(0,k_o)$ contrarily to $\alpha_j$.
Therefore, a horizontal array has reduced performance compared to a vertical array with the same length,
because the number of singular vectors that can be extracted for a given signal-to-noise ratio is reduced.

\begin{figure}
\begin{center}
\begin{tabular}{cc}
\includegraphics[width=6.cm]{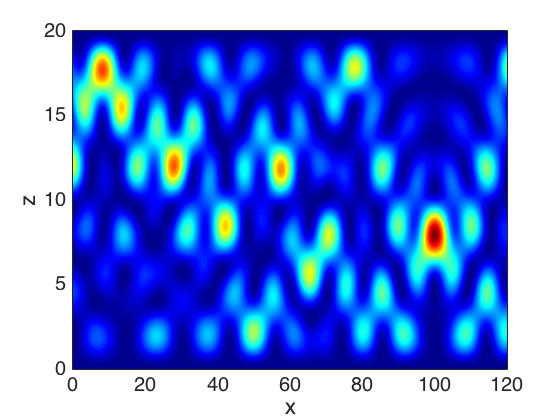}&
\includegraphics[width=6.cm]{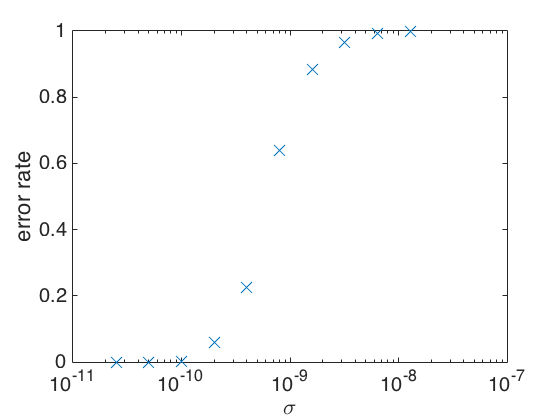}\\
$\sigma=0$ & 
localization error rate
\\
\includegraphics[width=6.cm]{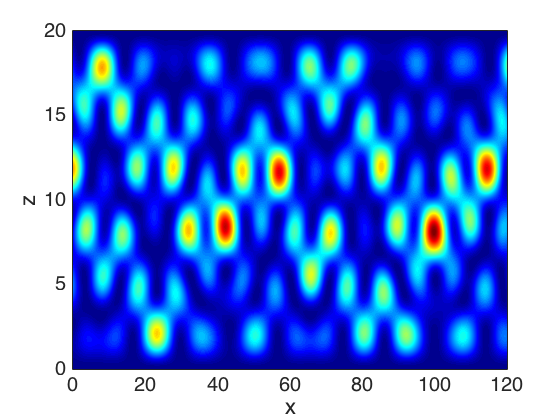}&
\includegraphics[width=6.cm]{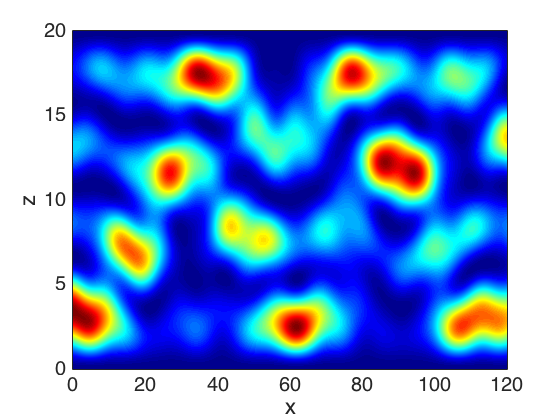} \\
$\sigma=10^{-9}$ & $\sigma=10^{-8}$
\end{tabular}
\end{center}
\caption{
{
Images and localization error rates obtained with a horizontal antenna array of $M=20$ receivers 
and total length $0.25 \simeq 0.04 \lambda_o$ and with different levels of noise~$\sigma$.
Here $\sigma_{\rm meas} = \sigma \| {\itbf p}\|_\infty$.
The source position is at $(x_o,z_o)=(100,7.7)$. 
The source can be localized if $\sigma \lesssim 10^{-10}$.
}
}
\label{fig:2}
\end{figure}

In figure \ref{fig:2},
an horizontal antenna array records the time-harmonic wave with different levels of additive noise. 
The array is localized at $z_{\rm a}=11$ and $x_k=0.25 (k-M/2)/M$, $k=1,\ldots,M$,
 with $M=20$.
Here the frequency is $\omega=1$, the velocity is $c_o=1$, $L=20$, the original source is at $(x_o,z_o)=(100,7.7)$.
The image is more sensitive to the noise than in the case of a vertical array, as predicted by the theory.

\subsection{Planar antenna array}
In this subsection we consider the situation in which the antenna array is planar and localized around $(x_{\rm a},z_{\rm a})$ 
and consists of $M$ receivers at $(x_k,z_k)$, $k=1,\ldots,M$.
The diameter of the planar array is denoted by $a$ (which means that the receivers lie within the square
$[x_{\rm a}-a,x_{\rm a}+a]\times [z_{\rm a}-a,z_{\rm a}+a]$).
The matrix  ${\bf B}$ has the form
\begin{equation}
{\bf B} = \big(  \phi_j(z_k)  \exp (-i \beta_j x_k) \big)_{1\leq k \leq M, 1\leq j \leq N}.
\end{equation}

As in the previous sections, by taking the regularization parameter $\epsilon$ proportional to the 
standard deviation of the meaurement noise, one minimizes the mean square estimation error. 
This means that we can estimate a limited number of pairs of singular values/vectors.
This number is the effective rank of the matrix ${\bf B}$, 
which depends on the antenna array and the waveguide geometry.
The goal of this subsection is to characterize this number and to show that it is not as large as one could have expected.

\subsubsection{Homogeneous waveguide}
Let us consider the case where the velocity is constant and the two boundary conditions are Dirichlet.
Then we have (\ref{eq:expressmodes}) and (\ref{eq:expressmodes2}) and
${\bf B} $ can be expanded as
\begin{equation}
\label{eq:expand:Bqqp}
{\bf B}=
\sum_{q,q'=0}^\infty 
{\itbf u}_{q,q'} {\itbf v}_{q,q'}^\dag  ,
\end{equation}
with
\begin{eqnarray}
{\itbf u}_{q,2q'+1} &=& \Big(  \frac{\sqrt{2}}{\sqrt{L}} \cos(\alpha_j z_{\rm a})  \beta_j^q\alpha_j^{2q'+1} \Big)_{j=1}^N,
\\
{\itbf u}_{q,2q'} &=&  \Big(  \frac{\sqrt{2}}{\sqrt{L}}\sin(\alpha_j z_{\rm a}) \beta_j^q\alpha_j^{2q'} \Big)_{j=1}^N,\\
{\itbf v}_{q,2q'+1} &=& \frac{i^q}{q!}\frac{(-1)^{q'+1}}{(2q'+1)!}
\Big( x_k^{q} \tilde{z}_k^{2q'+1}\Big)_{k=1}^M,
\\
{\itbf v}_{q,2q'} &=& \frac{i^q}{q!}\frac{(-1)^{q'}}{(2q')!}
\Big( x_k^{q} \tilde{z}_k^{2q'}\Big)_{k=1}^M .
\end{eqnarray}

\begin{proposition}
\label{lem:3}
 The vectors $({\itbf v}_{q,q'})_{q+q'\leq Q-1}$ are linearly independent
for arbitrary positions  $(x_k,\tilde{z}_k)_{k=1}^M$.
\end{proposition}
``Arbitrary" positions means outside a set of values of  $(x_k,\tilde{z}_k)_{k=1}^M$ in $[-a,a]^{2M}$ of Lebesgue measure zero.
For instance, if $(x_k,\tilde{z}_k)_{k=1}^M$ are sampled independently and randomly with the uniform distribution in $[-a,a]^2$,
then almost every realizations are ``arbitrary".  \\
{\it Proof.}
If the vectors $({\itbf v}_{q,q'})_{q+q'\leq Q-1}$ are linearly dependent, then the linear system
 $\sum_{q+q'\leq Q-1} \lambda_{q,q'} {\itbf v}_{q,q'}={\bf 0}$ has a non trivial solution $\boldsymbol{\lambda}= (\lambda_{q,q'})_{q+q'\leq Q-1}$.
 This system reads as ${\bf M} \boldsymbol{\lambda}= {\bf 0}$, where ${\bf M}$ is 
 a $M \times Q(Q+1)/2$-matrix involving the coefficients $\frac{i^q}{q!} \frac{(-1)^{[(q'+1)/2]}}{q'!} x_k^q \tilde{z}_k^{q'}$.
If $M \geq Q(Q+1)/2$ we can extract the first $Q(Q+1)/2$ lines of this matrix and we can claim that the determinant of
the resulting matrix must be zero. This determinant is a multivariate polynomial in $(x_k,\tilde{z}_k)_{k=1}^M$.
However, the set of points in $\mathbb{R}^{2M}$ 
in which a nonzero multivariate polynomial vanishes has zero Lebesgue measure \cite{caron}. \qed

If $k_oa$ is small ($a$ is the diameter of the planar antenna array), 
then
\begin{equation}
{\bf B}=
\sum_{q+q' \leq Q-1}
{\itbf u}_{q,q'} {\itbf v}_{q,q'}^\dag
+O\Big( \frac{ (k_oa)^{Q}}{Q!} \Big)  .
\end{equation}
The $\epsilon$-regularization used when the noise level has relative standard deviation $\epsilon$
prevents from exploiting the singular vectors whose singular values are smaller than 
$\epsilon$. 
Then it seems that the effective rank of ${\bf B}$ could be $(Q+1)Q/2$, i.e. the number of pairs $(q,q')$ such that $q+q'\leq Q-1$
where $Q$ is such that $\frac{ (k_oa)^{Q}}{Q!}\simeq \epsilon$.
Unfortunately, this is over-optimistic. 
Indeed, the vectors $({\itbf v}_{q,q'} )_{q+q'\leq Q-1}$ are
typically  linearly independent
by Proposition \ref{lem:3},  but the vectors $({\itbf u}_{q,q'} )_{q+q'\leq Q-1}$
are not.
Indeed, 
by the dispersion relation we have
\begin{equation}
\alpha_j^2+\beta_j^2=k_o^2,
\end{equation}
 which implies by Proposition \ref{lem:rank2dim}
\begin{equation}
{\rm Span}\Big( ({\itbf u}_{q,q'} )_{q+q'\leq Q-1} \Big)
={\rm Span}\Big( ({\itbf u}_{0,q'} )_{q' \leq Q-1}
\cup ({\itbf u}_{1,q'} )_{q'\leq Q-2} \Big)  ,
\end{equation}
and therefore the effective rank of ${\bf B}$ is only $2Q-1$.
However, $2Q-1$ is approximately twice as large as the effective
rank obtained for a horizontal array or a vertical one.
It is therefore much more favorable to use this type of antenna array.
But it is not as favorable as we could have anticipated.

We can therefore claim that ${\bf B} $ has an effective rank $2Q-1$
where $Q$ is such that $\frac{ (k_oa)^{Q}}{Q!}\simeq \epsilon$.
For instance, if $\epsilon=10^{-4}$ and $k_oa=0.125$, then $Q\simeq 3$ and the rank should be $5$.
If we compare with the singular values of the matrix ${\bf B}$ when $k_o =1 $, $L=20$, 
and $(x_k,z_k)_{k=1}^{M}$ is a Latin Hypercube Sampling (LHS) design (a type of quasi Monte Carlo sampling \cite{owen}) 
of $M=20$ points centered at $(0,z_{\rm a})$, $z_{\rm a}=11$, with size $0.25$,
then we find that the first singular values of ${\bf B}$ are $\sigma(1)\simeq 2.4$, $\sigma(2)\simeq 0.12$,
$\sigma(3)=0.03$, $\sigma(4)=2 \, 10^{-3}$, $\sigma(5)=1 \, 10^{-4}$, $\sigma(6)=5 \, 10^{-6}$,
so that indeed its effective rank is $5$.
 We observe that the singular values decay much slower than in the case of a vertical or horizontal array,
 which makes it possible to get an approximation of the inverse of  the matrix ${\bf B}$ even in the presence of moderate noise.

\begin{figure}
\begin{center}
\begin{tabular}{cc}
\includegraphics[width=6.cm]{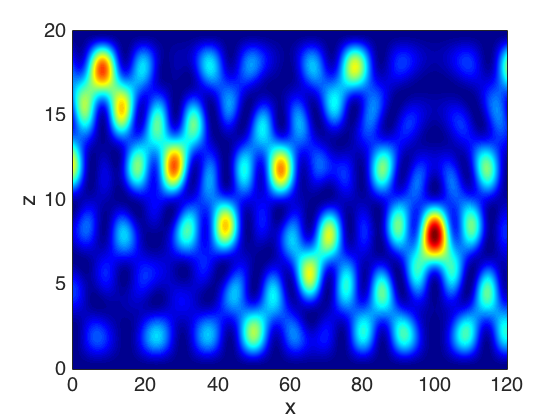} &
\includegraphics[width=6.cm]{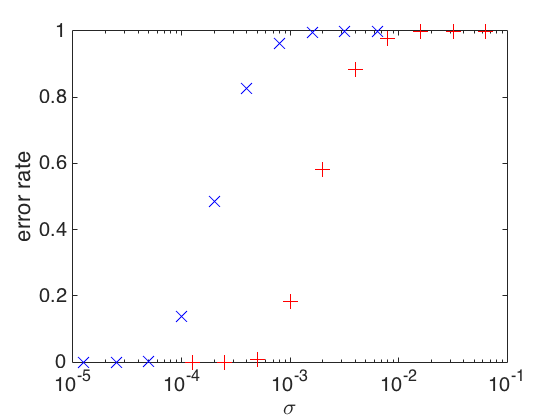} \\
$\sigma=0$ &
localization error rate
\\
\includegraphics[width=6.cm]{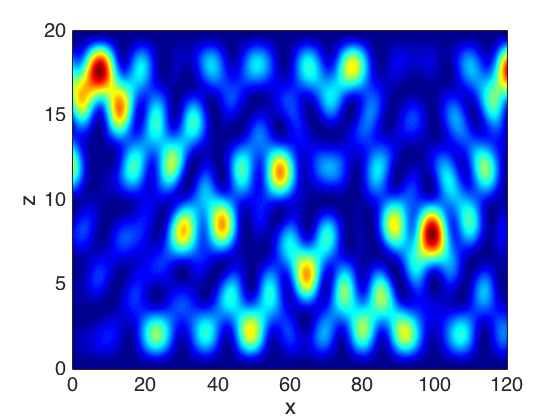}&
\includegraphics[width=6.cm]{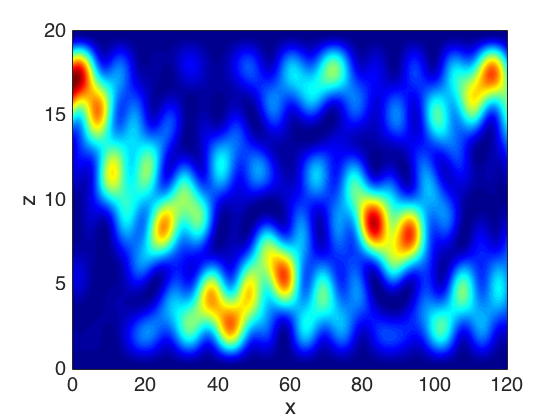}\\
$\sigma=10^{-4}$
&
$\sigma=10^{-3}$
\end{tabular}
\end{center}
\caption{
{
Images and localization error rates obtained with a planar array of $M=20$ receivers 
with side length $0.25\simeq 0.04 \lambda_o$ and with different levels of noise~$\sigma$.
Here $\sigma_{\rm meas} = \sigma \| {\itbf p}\|_\infty$.
The source position is at $(x_o,z_o)=(100,7.7)$ and the frequency is $\omega=1$. 
The localization error rate is plotted as a function of the noise level $\sigma$ for $M=20$ (blue $\times$) and for $M=1000$ (red $+$).
The source can be localized if $\sigma \lesssim 10^{-4}$ (for $M=20$) and if $\sigma \lesssim 10^{-3}$ (for $M=1000$).
}
}
\label{fig:3}
\end{figure}

\begin{figure}
\begin{center}
\begin{tabular}{cc}
\includegraphics[width=6.cm]{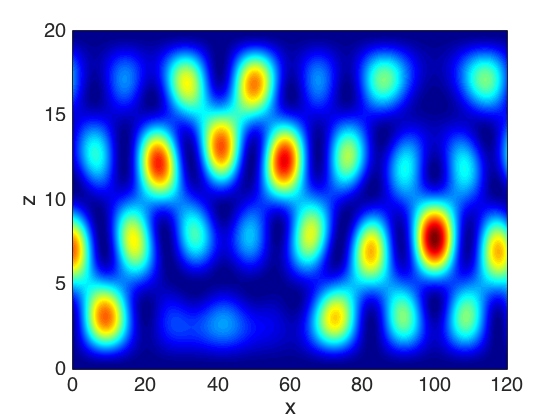} &
\includegraphics[width=6.cm]{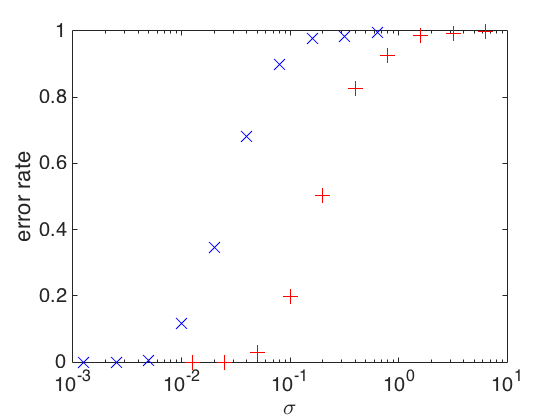} \\
$\sigma=0$ & localization error rate \\
\includegraphics[width=6.cm]{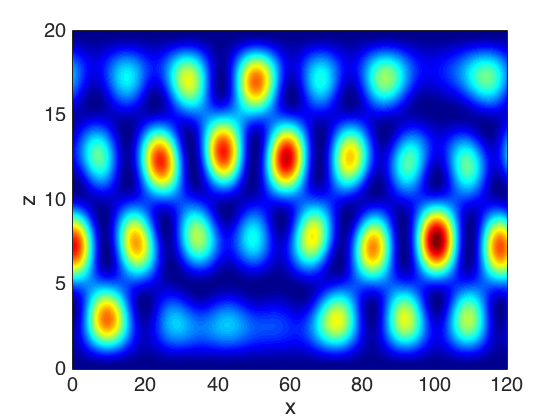}&
\includegraphics[width=6.cm]{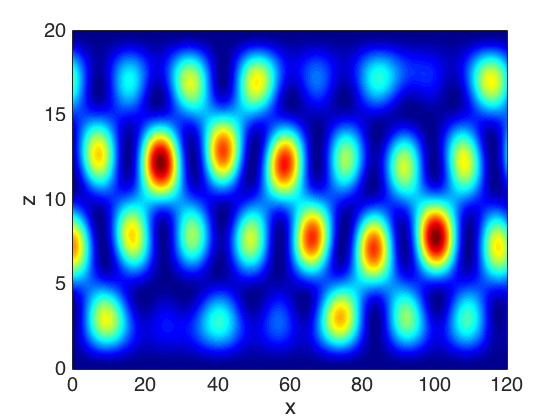}\\
$\sigma=10^{-2}$
&
$\sigma=10^{-1}$
\end{tabular}
\end{center}
\caption{
{
Images and localization error rates obtained with a planar array of $M=20$ receivers 
with side length $0.25\simeq 0.03 \lambda_o$ and with different levels of noise~$\sigma$. 
Here $\sigma_{\rm meas} = \sigma \| {\itbf p}\|_\infty$.
The source position is at $(x_o,z_o)=(100,7.7)$ and the frequency is $\omega=0.7$. 
The localization error rate is plotted as a function of the noise level $\sigma$ for $M=20$ (blue $\times$) and for $M=1000$ (red $+$).
The source can be localized if $\sigma \lesssim 10^{-2}$ (for $M=20$) and if $\sigma \lesssim 10^{-1}$ (for $M=1000$).
}
}
\label{fig:6}
\end{figure}

In figure \ref{fig:3},
a planar antenna array records the time-harmonic wave. It is centered at $(0,z_{\rm a})=(0,11)$.
It contains $M=20$ receivers distributed as a LHS design with size $a=0.25$.
Here the frequency is $\omega=1$, the velocity is $c_o=1$, $L=20$, the original source is at $(x_o,z_o)=(100,7.7)$.
The planar array can localize the source with a higher level of noise
compared to the linear, horizontal or vertical, array.
{We can observe that the source can be localized if $\sigma \lesssim 10^{-4}$.
If the number of receivers is multiplied by $K$ (for instance, $K=50$ as in figure \ref{fig:3} where the configurations with $M=20$ receivers and with $M=1000$ receivers are compared), then we gain a factor $\sqrt{K}$ 
in the critical noise level below which we can estimate the source position with accuracy and with high probability.
This gain can be observed in the figure plotting the localization error rates
by comparing the red ($M=20$) and blue ($M=1000$) crosses.}
Finally, if the frequency is $\omega=0.7$, then there is only $N=4$ guided modes
and it is possible to get the source position with a few percents of additive noise (see figure \ref{fig:6}).

\subsubsection{Heterogeneous waveguide}
The case of a homogeneous waveguide is special.
We can wonder whether the fact that  the effective rank of the matrix ${\bf B}$ is found to be $2Q-1$ while
we could have expected $Q(Q+1)/2$ is a particular feature of this waveguide or whether
it holds true for a general waveguide.
In fact, it turns out that this is a general feature that happens for any waveguide.
Indeed, in the general case, the matrix ${\bf B}$ can be expanded as
\begin{equation}
B_{kj}=
\sum_{q,q'=0}^\infty \phi_j^{(q')}(z_{\rm a}) 
 \frac{ \tilde{z}_k^{q'}}{q'!}
\frac{(-i\beta_j x_k)^{q}}{q!}  ,
\end{equation}
that is to say as  (\ref{eq:expand:Bqqp}) with
\begin{eqnarray}
{\itbf u}_{q,q'} &=& 
\Big( \beta_j^q \phi_j^{(q')}(z_{\rm a})  \Big)_{j=1}^N,
\\
{\itbf v}_{q,q'} &=& \frac{i^{q}}{q! q'!}
\Big( x_k^{q} \tilde{z}_k^{q'}\Big)_{k=1}^M .
\end{eqnarray}
If $k_oa$ is small ($a$ is the diameter of the planar array), 
then
\begin{equation}
{\bf B}=  
\sum_{q+q' \leq Q-1} 
{\itbf u}_{q,q'} {\itbf v}_{q,q'}^\dag
+O\Big( \frac{ (k_oa)^{Q}}{Q!} \Big) ,
\end{equation}
so that the effective rank of the matrix could be $Q(Q+1)/2$, provided the vectors ${\itbf v}_{q,q'}$ and 
the vectors ${\itbf u}_{q,q'}$, for $q+q'\leq Q-1$, are linearly independent (by assuming that $M,N$ are larger than $Q(Q+1)/2$).
For arbitrary positions $x_k, \tilde{z}_k$, the vectors ${\itbf v}_{q,q'}$ are linearly independent.
However, the vectors ${\itbf u}_{q,q'} $ are not independent,
as shown by the following proposition.
\begin{proposition}
\label{lem:rank2dim}
If $c$ is smooth at $z_{\rm a}$, then 
\begin{equation}
{\rm Span}\Big( ({\itbf u}_{q,q'} )_{q+q'\leq Q-1} \Big)
={\rm Span}\Big( ({\itbf u}_{q,0} )_{q \leq Q-1}
\cup ({\itbf u}_{q,1} )_{q\leq Q-2} \Big)  .
\end{equation}
\end{proposition}
\noindent
{\it Proof.}
Eq.~(\ref{eq:vp}) implies the following relations for the derivatives of the mode profiles:
\begin{eqnarray*}
\phi_j^{(2)}(z_{\rm a}) &=& \, \big(\beta_j^2 - \frac{\omega^2}{c^2(z_{\rm a})} \big) \phi_j(z_{\rm a}) ,\\
\phi_j^{(3)}(z_{\rm a}) &=& \, \frac{2 \omega^2 c^{(1)}(z_{\rm a})}{c^3(z_{\rm a})} 
\phi_j(z_{\rm a})
+
\big(\beta_j^2- \frac{\omega^2}{c^2(z_{\rm a})}\big)  \phi_j^{(1)}(z_{\rm a}) ,\\
\phi_j^{(4)}(z_{\rm a}) &=& \, \Big(  \big( \beta_j^2 -\frac{\omega^2}{c(z_{\rm a})^2}\big)^2
+
 \omega^2 \Big( \frac{2c^{(2)}(z_{\rm a})}{c^3(z_{\rm a})} - \frac{6c^{(1)}(z_{\rm a})^2}{c^4(z_{\rm a})} \Big) \Big)
\phi_j(z_{\rm a})\\
&&+ \frac{4 \omega^2c^{(1)}(z_{\rm a})}{c^3(z_{\rm a})}  \phi_j^{(1)}(z_{\rm a}),
\end{eqnarray*}
and more generally we can establish by a recursive argument that
\begin{equation}
\phi_j^{(q)}(z_{\rm a}) = P_q(\beta_j^2) \phi_j(z_{\rm a}) +Q_q(\beta_j^2) \phi_j^{(1)}(z_{\rm a})  ,
\end{equation}
where $P_q$ and $Q_q$ are polynomials of degree ${\rm deg}(P_q) \leq [q/2]$ and ${\rm deg}(Q_q)\leq [(q-1)/2]$, whose coefficients depend on $\omega$ and
on derivatives $c^{(k)}(z_{\rm a})$, but not explicitly on $j$:
\begin{eqnarray*}
P_{q+1}(\beta^2) &=& [\partial_{z_{\rm a}}P_q] (\beta^2) +Q_q(\beta^2) \big(\beta^2 -\frac{\omega^2}{c^2(z_{\rm a})}\big) ,\\
Q_{q+1}(\beta^2) &=& P_q (\beta^2) +[\partial_{z_{\rm a}}Q_q] (\beta^2) .
\end{eqnarray*}
Consequently
$$
 \beta_j^q \phi_j^{(q')}(z_{\rm a}) = \big[  \beta_j^q P_{q'}(\beta_j^2) \big] \phi_j(z_{\rm a}) + \big[ \beta_j^q Q_{q'}(\beta_j^2) \big] \phi_j^{(1)}(z_{\rm a}) .
$$
This shows that ${\itbf u}_{q,q'}$ is a linear combination of  $({\itbf u}_{q+2q'',0})_{q''=0}^{[q'/2]}$ and $({\itbf u}_{q+2q'',1})_{q''=0}^{[(q'-1)/2]}$,
and it is therefore a linear combination of $({\itbf u}_{q'',0})_{q''=0}^{q+q'}$ and $({\itbf u}_{q'',1})_{q''=0}^{q+q'-1}$.
\qed

Therefore the effective rank of ${\bf B}$ is only $2Q-1$.
Again, the effective rank $2Q-1$ is twice as large as the effective
rank obtained for a horizontal array or a vertical one.
It is therefore much more favorable to use this type of antenna array.
Unfortunately, it is not possible to reach the rank $Q(Q+1)/2$ that would have been even more favorable.

\begin{figure}
\begin{center}
\begin{tabular}{cc}
\includegraphics[width=6.cm]{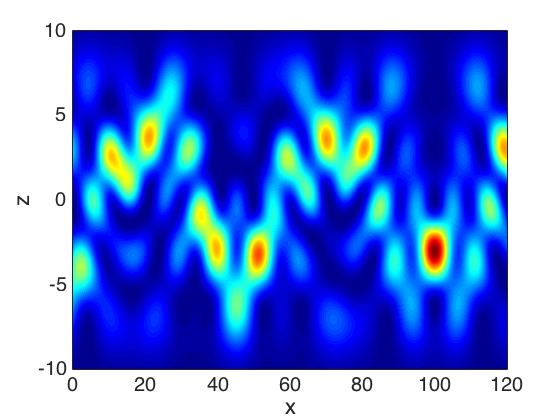} &
\includegraphics[width=6.cm]{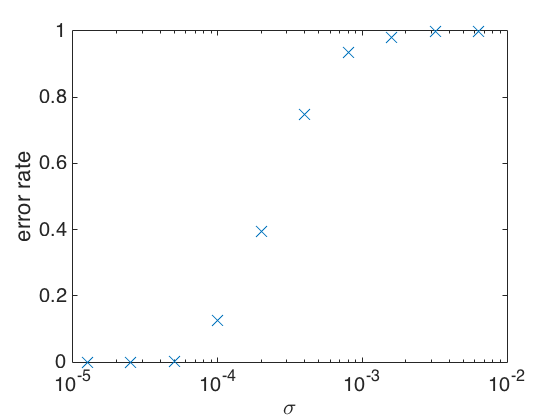} \\
$\sigma=0$ &  localization error rate \\
\includegraphics[width=6.cm]{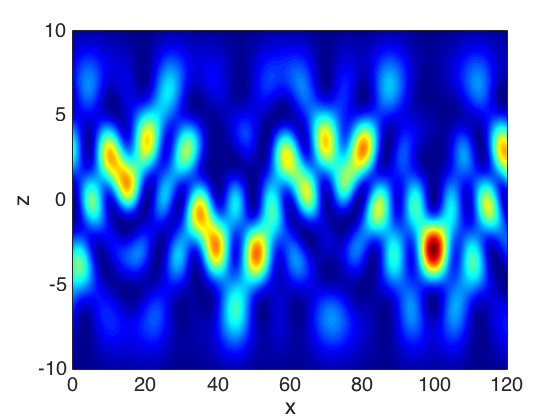}&
\includegraphics[width=6.cm]{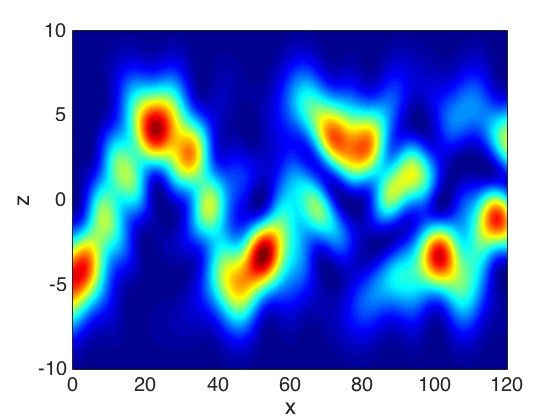}\\
$\sigma=10^{-4}$
&
$\sigma=10^{-3}$
\end{tabular}
\end{center}
\caption{
{
Images and localization error rates obtained with a planar array of $M=20$ receivers 
with side length $0.25\simeq 0.04 \lambda_o$ and with different levels of noise~$\sigma$.
Here $\sigma_{\rm meas} = \sigma \| {\itbf p}\|_\infty$.
The waveguide is parabolic. 
The source position is at $(x_o,z_o)=(100,-3)$. 
The source can be localized if $\sigma \lesssim 10^{-4}$.
}
}
\label{fig:5}
\end{figure}

\begin{example} 
Let us consider the case of an ideal parabolic waveguide, with unbounded transverse domain and 
a transverse velocity profile of the form
$$
\frac{1}{c^2(z)} = \frac{1}{c_o^2} \Big(1 - \frac{z^2}{L^2}\Big) .
$$
Denoting $k_o=\omega/c_o$, 
the eigenmodes have the form
$$
\phi_j(z) = (k_o/L)^{1/4} f_j\big( (k_o/L)^{1/2}z\big) ,
$$
where the $f_j$, $j \geq 0$, are the Gauss-Hermite functions:
$$
f_j(s) = \frac{1}{\sqrt{2^j \sqrt{\pi}j!}}H_j(s)\exp(-s^2/2),
$$
that satisfy $f_j''(s) -s^2 f_j (s)=-(2j+1) f_j(s)$. 
There are $N+1$ guided modes, with
$$
N = \left[ \frac{k_oL-1}{2}\right] .
$$
For $j=0,\ldots,N$, the modal wavenumber of the $j$th guided mode is
$$
\beta_j =  \sqrt{k_o^2 -(2j+1) k_o/L} .
$$
In figure \ref{fig:5},
a planar antenna array records the time-harmonic wave. It is centered at $(0,z_{\rm a})=(0,2)$.
It contains $M=20$ receivers distributed as a LHS design with size $a=0.25$.
Here the frequency is $\omega=1$, the waveguide is parabolic with $c_o=1$ and $L=10$, the original source is at $(x_o,z_o)=(100,-3)$.
The result is very similar to the case of a homogeneous waveguide.
\end{example}

\section{Conclusion}
In this paper we have considered the source imaging problem in a two-dimensional waveguide.
We have first addressed the case of dense antenna arrays 
in the high-frequency regime and compared the performances of
vertical and horizontal antenna arrays. The overall result is that the length of a horizontal antenna array 
should be twice as long as the one of a vertical array to present similar performance.
We have focused our attention to the low-frequency regime,
when the number of guided modes is small and the diameter of the sensor array is smaller than the wavelength.
The principle of source localization is 1) to estimate the guided mode amplitudes from the recorded data by resolution of an 
appropriate regularized inverse problem and 2) to backpropagate the contributions of the guided mode amplitudes that have been 
estimated correctly.
The main findings of this paper are the following ones: \\
i) Source localization is possible even with very small antenna arrays provided the signal-to-noise ratio (SNR) of the data is high.\\
ii) Vertical linear antenna arrays have better performance than horizontal linear arrays (for a given diameter) but both require extremely high SNR.\\
iii) The use of planar antenna arrays makes it possible to get an estimate of the source position when the SNR is moderately high.
The gain in performance and stability compared to linear (horizontal or vertical) antenna arrays is significant.\\
iv) There is a fundamental limitation that prevents from reaching an even better performance
and that is related to the wave equation and its dispersion relation. 
It is one situation where a PDE-constrained inverse problem (an inverse problem constrained by a partial differential equation)
shows poor results because the acquired data set is in fact highly redundant.

\section*{Acknowledgements}
This work was partly supported by Direction G\'en\'erale de l'Armement (DGA) Naval Systems and by
ANR under Grant No.~ANR-19-CE46-0007 (project ICCI).


\section*{References}

\begin{thebibliography}{9}

\bibitem{ammarigarniersolna}
{
H. Ammari, J. Garnier, and K. S\o lna, 
A statistical approach to target detection and localization in the presence of noise, 
Waves Random Complex Media {\bf 22} (2012), 40--65.
}

\bibitem{arens11}
T. Arens, D. Gintides, and A. Lechleiter,
Direct and inverse medium scattering in a three-dimensional homogeneous planar waveguide,
SIAM J. Appl. Math. {\bf 71} (2011), 753--772.

\bibitem{borcea10}
L. Borcea, T. Callaghan, J. Garnier, and G. Papanicolaou, 
A universal filter for enhanced imaging with small arrays, 
Inverse Problems {\bf 26} (2010), 01506.

\bibitem{bgt15}
L. Borcea, J. Garnier, and C. Tsogka,
A quantitative study of source imaging in random waveguides,
Commun. Math. Sci.  {\bf 13} (2015), 749--776.

\bibitem{bourgeois08}
L. Bourgeois and E. Lun\'eville,
The linear sampling method in a waveguide: a modal formulation,
Inverse problems {\bf 24} (2008), 015018.

\bibitem{buchanan}
J. L. Buchanan, R. P. Gilbert, A. Wirgin, and Y. Xu, 
Marine Acoustics: Direct and Inverse Problems, SIAM, Philadelphia, 2004.

\bibitem{caron}
R. Caron and T. Traynor, The zero set of a polynomial, WSMR Report 05-02, 2005.

\bibitem{dediu06}
S. Dediu and J. R. McLaughlin,
Recovering inhomogeneities in a waveguide using eigensystem decomposition,
Inverse Problems {\bf 22} (2006), 1227--1246.

\bibitem{engl}
H. W. Engl, M. Hanke, and A. Neubauer,
Regularization of Inverse Problems, 
volume 375 of Mathematics and its Applications, 
Kluwer, Dordrecht, 1996.

\bibitem{fasino}
D. Fasino,
Spectral properties of Toeplitz-plus-Hankel matrices,
Calcolo {\bf 33} (1996), 87--98.

\bibitem{book}
J.-P. Fouque, J. Garnier, G. Papanicolaou, and K. S\o lna, 
Wave Propagation and Time Reversal in Randomly Layered Media, Springer, New York, 2007.

\bibitem{gar-pap-06}
{\rm J. Garnier and G. Papanicolaou}, 
{Pulse propagation and time reversal in random waveguides},
SIAM J. Appl. Math. {\bf  67} (2007), 1718--1739.

\bibitem{garnierpapa}
{
{\rm J. Garnier and G. Papanicolaou}, 
Passive Imaging with Ambient Noise,
Cambridge University Press, Cambridge, 2016.
}

\bibitem{grenander}
U. Grenander and G. Szeg\"o, 
Toeplitz Forms and Their Applications, 
Chelsea, New York, 1984.

\bibitem{hastie}
T. Hastie, R. Tibshirani, and J. Friedman,
 The Elements of Statistical Learning,
 Springer, New York, 2009.
 
\bibitem{hodgkiss99}
W. S. Hodgkiss, H. C. Song, W. A. Kuperman, T. Akal, C. Ferla, and D. R. Jackson,
A long-range and variable focus phase-conjugation experiment in shallow water,
J. Acoust. Soc. Amer.  {\bf 105} (1999), 1597--1604.

\bibitem{jensen11}
F. B. Jensen, W. A. Kuperman, M. B. Porter, and H. Schmidt, 
Computational Ocean Acoustics,
 Springer, New York, 2011.

\bibitem{kuperman02}
W. A. Kuperman and D. Jackson,
Ocean acoustics, matched-field processing and phase conjugation,
 Topics in Applied Physics, pages 43--97, Springer, Berlin, 2002.

\bibitem{monk12}
P. Monk and V. Selgas,
Sampling type methods for an inverse waveguide problem,
Inverse Problems and Imaging {\bf 6} (2012), 709--747.

\bibitem{mordant}
N. Mordant, C. Prada, and M. Fink, 
Highly resolved detection in a waveguide using the D.O.R.T. method, 
J. Acoust. Soc. Amer. {\bf 105} (1999), 2634--2642.

\bibitem{owen}
A. B. Owen, 
Orthogonal arrays for computer experiments, integration and visualization,
Statistica Sinica {\bf 2} (1992), 439--452.

\bibitem{pincon}
B. Pin\c{c}on and K. Ramdani, 
Selective focusing on small scatterers in acoustic waveguides using time reversal mirrors, 
Inverse Problems {\bf 23} (2007), 1--25.

\bibitem{prada07}
C. Prada, J. de Rosny, D. Clorennec, J.-G. Minonzio, A. Aubry, M. Fink, L. Berniere, P. Billand, S. Hibral, and T. Folegot, 
Experimental detection and focusing in shallow water by decomposition of the time reversal operator, 
J. Acoust. Soc. Amer. {\bf 122} (2007), 761--768.

\bibitem{scherzer}
O. Scherzer, M. Grasmair, H. Grossauer, M. Haltmeier, and F. Lenzen,
Variational Methods in Imaging, 
volume 167 of Applied Mathematical Sciences, Springer, New York, 2009.

\bibitem{sengupta}
I. SenGupta, B. Sun, W. Jiang, G. Chen, and M. C. Mariani,
Concentration problems for bandpass filters in communication theory over disjoint frequency intervals and numerical solutions,
J. Fourier Anal. Appl. {\bf 18} (2012), 182--210.

\bibitem{slepian32}
D. Slepian,
Prolate spheroidal wave functions, Fourier analysis, and uncertainty - V: The discrete case,
Bell System Technical Journal {\bf 57} (1978), 1371--1430.

\bibitem{tsogka13}
C. Tsogka, D. A. Mitsoudis, and S. Papadimitropoulos,
Selective imaging of extended reflectors in two-dimensional waveguides,
SIAM J. Imaging Sci. {\bf 6} (2013), 2714--2739.

\bibitem{tsogka16}
C. Tsogka, D. A. Mitsoudis, and S. Papadimitropoulos, 
Partial-aperture array imaging in acoustic waveguides, 
Inverse Problems {\bf 32} (2016), 125011.

\bibitem{tsogka18}
C. Tsogka, D. A. Mitsoudis, and S. Papadimitropoulos, 
Imaging extended reflectors in a terminating waveguide, 
SIAM J. Imaging Sci. {\bf 11} (2018), 1680--1716. 


\end{thebibliography}
\end{document}